\documentclass[11pt,a4paper,twoside]{article}
\usepackage{color}
\usepackage{bm, amsmath, amssymb, amsthm} 

\def\calf{{\cal F}}
\def\<{\langle}
\def\>{\rangle}
\def\eps{\varepsilon}
\def\RR{\mathbb{R}}

\newcommand\const{\operatorname{const}}
\newcommand\tr{\operatorname{Tr}}
\newcommand\Div{\operatorname{div}}

\def\Ric{\operatorname{Ric}}
\def\vol{\operatorname{vol}}

\def\eq{\hspace*{-1.5mm}&=&\hspace*{-1.5mm}}
\def\plus{\hspace*{-1.5mm}&+&\hspace*{-1.5mm}}
\def\minus{\hspace*{-1.5mm}&-&\hspace*{-1.5mm}}
\def\dt{\partial_t}

\newtheorem{corollary}{Corollary}
\newtheorem{definition}{Definition}

\newtheorem{remark}{Remark}
\newtheorem{lemma}{Lemma}
\newtheorem{proposition}{Proposition}
\newtheorem{theorem}{Theorem}

\author{Vladimir Rovenski\footnote{Mathematical Department, University of Haifa, Mount Carmel, 31905 Haifa,  Israel
       \newline e-mail: {\tt rovenski@math.haifa.ac.il}
       }
        \ and \
        Pawe\l\ Walczak\footnote{Katedra Geometrii,
        Uniwersytet \L\'{o}dzki, ul. Banacha 22,
             90-238  \L\'{o}d\'{z}, Poland
        \newline e-mail: {\tt pawelwal@math.uni.lodz.pl}}
}

\title{Integral formulae\\ for codimension-one foliated Finsler manifolds}

\begin{document}

\date{}

\maketitle


\begin{abstract}
We study extrinsic geometry of a codimension-one foliation $\calf$ of a
Finsler space $(M,F)$, in particular, of a Randers space $(M,\alpha+\beta)$.
Using a unit vector field $\nu$ orthogonal (in the Finsler sense) to the leaves of $\calf$,
we define a new Riemannian metric $g$ on $M$, which for Randers case depends nicely on $(\alpha,\beta)$.
For that $g$ we derive several geometric invariants of $\calf$ (e.g. the Riemann curvature and the shape operator)
in terms of $F$; then under natural assumptions on $\beta$ which simplify derivations,
we express them in terms of invariants arising from $\alpha$ and~$\beta$.
Using our approach of \cite{rw2}, we produce the integral formulae for $\calf$ of closed $(M, F)$ and $(M, \alpha+\beta)$,
which relate integrals of mean curvatures with those involving algebraic invariants obtained from
the shape operator of $\calf$ and the Riemann curvature in the direction $\nu$.
They generalize formulae
by Brito-Langevin-Rosenberg (that total mean curvatures of any order for a
foliated closed Riemannian space of constant curvature don't depend on a choice of~$\calf$).


\end{abstract}

\section*{Introduction}

Two recent decades brought increasing interest in Finsler geometry (see \cite{bcs,cs,sh2}
and the bibliographies therein), in particular, in extrinsic geometry of hypersurfaces
of Finsler manifolds (see the items above and, for example, \cite{sh1}).
Among all the Finsler structures, Randers metrics (introduced in \cite{ra}
and being the closest relatives of Riemannian ones) play an important role.

Extrinsic geometry of foliated Riemannian manifolds is also of definite interest
since some time (see \cite{ro,rw1} and, again, the bibliographies therein).
Among other topics of interest, one can find a number of papers devoted to so called
{\it integral formulae} (see surveys in  \cite{rw1,arw2014}), which provide obstructions for existence of foliations
(or compact leaves of them) with given geometric properties.
A series of integral formulae has been provided in \cite{rw2}. They include the formulae in
\cite{re} that the~total mean curvature of the leaves is zero, and generalize the formulae in \cite{blr}, which show that
total mean curvatures (of arbitrary order $k$) for codimension-one foliations
on a closed $(m+1)$-dimensional manifold of constant sectional curvature $K$ depend only
on $K$, $k$, $m$ and the volume of the manifold, not on a foliation.
One of such formulae was used in \cite{lw} to prove that codimension-one foliations
of a closed Riemannian manifold of negative Ricci curvature are far (in a sense defined there) from being~umbilical.

In this paper we study extrinsic geometry of a codimension-one transversely oriented foliation $\calf$ of a closed
Finsler space $(M,F)$, in particular, of a Randers space $(M, \alpha+\beta)$,
$\alpha$ being the norm of a Riemannian structure $a$ and $\beta$ a 1-form of $\alpha$-norm smaller
than $1$ everywhere on~$M$. Using a unit normal $\nu$ (in the Finsler sense)
to the leaves of $\calf$ we define a new Riemannian structure $g$ on $M$, which
in Randers case depends nicely on $\alpha$ and $\beta$.
 For that $g$, we derive several geometric invariants of $\calf$ (e.g. the Riemann curvature and the shape operator)
in terms of $F$;
under natural assumptions on $\beta$ which simplify derivations,
we express them in terms of corresponding invariants arising from $\alpha$ and some quantities related to $\beta$.
Then, using the approach of \cite{rw2}, we produce the integral formulae for $\calf$ on $(M, F)$ and $(M, \alpha+\beta)$;
some of them generalize the formulae in \cite{blr}.

Our formulae relate integrals of $\sigma_i$'s with those involving algebraic invariants (see Appendix)
obtained from $A_p$ ($p\in M$) -- the shape operator of a foliation $\calf$,
$R_p$ -- the Riemann curvature in the direction $\nu$ normal to $\calf$,
and their products of the form $(R_p)^j A_p$, $j=1,2,\ldots$
In fact, we get a bit more: we produce an infinite sequence of such formulae
for a smooth unit vector field $\nu$ on $M$ involving these algebraic invariants.
To simplify calculations, we work on locally symmetric ($\nabla R=0$ with respect to $g$) Finsler manifolds,
where our approach can be applied with the full force (Section~\ref{sec:intform}).
We show that our formulae reduce to these in \cite{blr} in the case of constant curvature
and to those in \cite{rw2} in the Riemannian case.
Using Finsler geometry of Randers spaces we produce  also (Section~\ref{sec:randers}) integral  formulae on
codimension-one foliated Riemannian manifolds which involve not only $A_p$'s and $R_p$'s but also an auxiliary 1-form $\beta$.

We discuss a number of particular cases and provide consequences of our new formulae.

\section{Preliminaries}

Recall Euler's Theorem: If a function $f$ on $\RR^{m+1}$ is smooth away from the origin of $\RR^{m+1}$
then the following two statements are equivalent:

-- $f$ is positively homogeneous of degree $r$, that is $f(\lambda\,y)=\lambda^r f(y)$ for all $\lambda>0$;

-- the radial derivative of $f$ is $r$ times $f$, namely, $f_{y^i}(y)\,y^i = r f(y)$.

\noindent
The obvious consequence of Euler's Theorem helps us to represent several formulae in what follows:

\begin{corollary}\label{L-hom-2}
If a smooth function $f$ on $\RR^{m+1}\setminus\{0\}$ obeys the 2-homogeneity condition $f(\lambda\,y)=\lambda^2 f(y)$ for $\lambda>0$ then
$f(y)= \frac12\,f_{y^i y^j}(y)\,y^i y^j $ for smooth functions $f_{y^i y^j}$ on $\RR^{m+1}\setminus\{0\}$.
\end{corollary}

\proof
By Euler's Theorem, $f_{y^i}(y)\,y^i = 2 f(y)$.
Since $f_{y^i}(\lambda\,y)=\lambda f_{y^i}(y)$, by Euler's Theorem, we have $f_{y^i}(y) = f_{y^i y^j}(y)y^j$.
\qed

\subsection{The Minkowski and Randers norms}
\label{sec:M-R-norms}

\begin{definition}[see \cite{sh2}]\rm
A \textit{Minkowski norm} on a vector space $\RR^{m+1}$ is a function $F:\RR^{m+1}\to[0,\infty)$ with the following properties
(of regularity, positive 1-homogeneity and strong convexity):

M$_1:$ $F\in C^\infty(\RR^{m+1}\setminus \{0\})$,\quad
M$_2:$ $F(\lambda\,y)=\lambda F(y)$ for all $\lambda>0$ and $y\in\RR^{m+1}$,

M$_3:$ For any $y\in \RR^{m+1}\setminus \{0\}$, the following symmetric bilinear form  is positive definite on $\RR^{m+1}:$
\begin{equation}\label{E-acsiom-M3}
 g_y(u,v)=\frac12\,\frac{\partial^2}{\partial s\,\partial t}\,\big[F^2(y+su+tv)\big]_{|\,s=t=0}\,.
\end{equation}
\end{definition}

\noindent
By (M$_2$), $g_{\lambda y}=g_{y}$ for all $\lambda>0$.
By~(M$_3$), $\{y\in \RR^{m+1}: F(y) \le 1\}$ is a strictly convex set.
Note that
\begin{equation}\label{E-g-F2}
 g_y(y,v)=\frac12\frac{\partial}{\partial t}\,\big[F^2(y+tv)\big]_{|\,t=0},\quad
 g_y(y,y)=F^2(y).
\end{equation}
One can check that $F(u+v)\le F(u) + F(v)$ (the triangle inequality)
and $F_{y^i}(y)\,u^i\le F(u)$  (the fundamental inequality) for all $y\in\RR^{m+1}\setminus\{0\}$ and $u, v\in\RR^{m+1}$.
By Corollary~\ref{L-hom-2}, we have
 $F^2(y)=g_{ij}(y)\,y^i y^j$, where $g_{ij}=\frac12\,[F^2]_{y^i y^j}=F F_{y^i y^j} +F_{y^i}F_{y^j}$
are smooth functions in $\RR^{m+1}\setminus\{0\}$ which, in general, cannot be extended continuously to all of $\RR^{m+1}$.
The following symmetric trilinear form $C$ for Minkowski norms is called the \textit{Cartan torsion}:
\begin{equation}\label{E-Ctorsion}
  C_y(u,v,w)=\frac12\,\frac{\partial}{\partial t}\,\big[g_{y+tw}(u,v)\big]_{|\,t=0}\quad
 {\rm where} \quad y\in\RR^{m+1}\setminus\{0\},\ u,v,w\in\RR^{m+1}\,.
\end{equation}
The homogeneity of $F$ implies the following:
\[
 C_y(u,v,w)=\frac14\,\frac{\partial^3}{\partial r\,\partial s\,\partial t}\,\big[F^2(y+ru+sv+tw)\big]_{|\,r=s=t=0},\quad
 C_{\lambda y}=\lambda^{-1}C_y\quad (\lambda>0).
\]
We have $C_y(y,\cdot\,,\cdot\,)=0$.
The \textit{mean Cartan torsion} is given by $I_{y}(u):=\tr C_{y}(\cdot\,,\cdot\,, u)$.
 Observe that
\[
 C_{ijk}:=C(\partial_{y^i},\partial_{y^j},\partial_{y^k})
 =\frac12\,\frac{\partial}{\partial y^k}\,g_{ij}=\frac14\,[F^2]_{y^i y^j y^k},\qquad
 I_{k}=g^{ij}C_{ijk}.
\]

Let $(b_i)$ be a basis for $\RR^{m+1}$ and $(\theta^i)$ the dual basis in $(\RR^{m+1})^*$.
The \textit{Busemann-Hausdorff volume form} is defined by ${\rm d}V_F=\sigma_F(x)\,\theta^1\wedge\dots\wedge\theta^{m+1}$, where $\sigma_F=\frac{\vol\mathbb{B}^{m+1}}{\vol B^{m+1}}$.
Here
$\mathbb{B}^{m+1}:=\{y\in\RR^{m+1}: \|y\|<1\}$ is a Euclidean unit ball, and $\vol B^{m+1}$ is the Euclidean volume
of a strongly convex subset $B^{m+1}:=\{y\in\RR^{m+1}: F(y^i b_i)<1\}$
(so that for the unit cubic $\mathcal{U}=[0,1]^{m+1}$, $\vol\,\mathcal{U}=1$).

The \textit{distortion} of $F$ is defined by
 $\tau(y)=\log({\sqrt{\det g_{ij}(y)}}/{\sigma_F})$.
It has the $0$-homogeneity property: $\tau(\lambda y)=\tau(y)\ (\lambda>0)$,
and $\tau=0$ for Riemannian spaces.

The \textit{angular form} is defined by $h_y(u,v)=g_y(u,v) - F(y)^{-2}g_y(y,u)\,g_y(y,v)$.
Observe that $h_y(u,u)\ge g_y(u,u) - F(y)^{-2}g_y(y,y)\,g_y(u,u)=0$ and equality holds if and only if $u||\,y$.


A~vector $n\in\RR^{m+1}$ is called a \textit{normal} to a hyperplane $W\subset\RR^{m+1}$ if
 $g_n(n,w)=0\ (w\in W)$.
There are exactly two normal directions to $W$, see \cite{sh2}, which are opposite when $F$ is \textit{reversible}
(i.e., $F(-y) = F(y)$ for all $y\in \RR^{m+1}$).

\begin{definition}\rm
 Let $a(\cdot\,,\cdot)= \<\cdot\,,\cdot\>$ be a scalar product and
 $\alpha(y)=\|y\|_\alpha=\sqrt{\<y,y\>}$ for $y\in \RR^{m+1}$ the corresponding Euclidean norm on $\RR^{m+1}$.
 If $\beta$ is a linear form on $\RR^{m+1}$ with $\|\beta\,\|_\alpha<1$ then the following function $F$
 is called the \textit{Randers norm}:
\begin{equation}\label{E-F000}
 F(y) = \alpha(y) + \beta(y) = \sqrt{\<y,y\>}+\beta(y).
\end{equation}
\end{definition}

For Randers norm \eqref{E-F000} on $\RR^{m+1}$, the bilinear form $g_y$ obeys, see \cite{sh2},
\begin{eqnarray}\label{E-F001}
\nonumber
 g_y(u,v)\eq \alpha^{-2}(y)(1+\beta(y))\,\<u,v\> +\beta(u)\,\beta(v) \\
 \minus \alpha^{-3}(y)\,\beta(y)\,\<y,u\>\,\<y,v\> +\alpha^{-1}(y)\,\big(\beta(u)\,\<y,v\> +\beta(v)\,\<y,u\>\big)\,,\\
\label{E-F001b}
 \det g_y \eq (F(y)/\alpha(y))^{m+2}\det a.
\end{eqnarray}
Let $N\in \RR^{m+1}$ be a unit normal to a hyperplane $W$ in $\RR^{m+1}$ with respect to $\<\cdot\,,\cdot\>$, i.e.,
\[
 \<N,w\>=0\quad (w\in W),\qquad
  \alpha(N)=\|N\|_\alpha=\sqrt{\<N,N\>}=1.
\]
Let $n$ be a vector $F$-normal to $W$, lying in the same half-space with $N$ and such that $\|n\|_\alpha=1$.
Set
\[
 g(u,v):=g_n(u,v),\quad u,v\in\RR^{m+1}.
\]
Then $g(n,n)=F^2(n)$, see \eqref{E-g-F2}, and  $F(n)=1+\beta(n)$.

The 'musical isomorphisms' $\sharp$ and $\flat$ will be used for rank one tensors
and symmetric rank 2 tensors on $(\RR^{m+1},a)$ and Riemannian manifolds.
For~example, if $\beta$ is a 1-form on $\RR^{m+1}$ and $v\in \RR^{m+1}$ then
$\<\beta^\sharp,u\>=\beta(u)$ and $v^\flat(u) =\<v,u\>$ for any $u\in\RR^{m+1}$.

\begin{lemma}\label{L-c-value}
 If the Randers norm obeys $\,\beta(N)=0$ $($i.e., $\beta^\sharp\in W)$ then
\begin{eqnarray}\label{E-c-value}
 n \eq c\,N-\beta^\sharp, \\
\label{E-c-value2}
 g(u,v) \eq c^2\big(\<u,v\> -\beta(u)\,\beta(v)\big),\quad u,v\in W\,,\\
\label{E-c-value2c}
 g(n,n) \eq c^4,\quad g(n,v) = 0,
\end{eqnarray}
where $c:=(1-\|\beta\,\|^2_\alpha)^{1/2}>0$.
 The vector $\nu=c^{-2}n$ is an $F$-unit normal to $W$.
\end{lemma}

\proof For arbitrary $\beta$ and $y=n$ and $\alpha(n)=1$, the formula \eqref{E-F001} reads
\begin{equation}\label{E-c-value0}
 g(u,v) = (1+\beta(n))\<u,v\> +\beta(u)\,\beta(v) - \beta(n)\,\<n,u\>\,\<n,v\> +\beta(u)\,\<n,v\> +\beta(v)\,\<n,u\>.
\end{equation}
Assuming $u=n$, from \eqref{E-c-value0} we find
\begin{equation}\label{E-gnv}
 g(n,v) =(1+\beta(n))\,\<n + \beta^\sharp,\,v\>.
\end{equation}
Note that
 $|\beta(n)|=|\<\beta^\sharp,n\>|\le\alpha(\beta^\sharp)\,\alpha(n)<1$; hence, $1+\beta(n)>0$.
We find from \eqref{E-gnv} with $v\in W$ that $n+\beta^\sharp=\hat c\,N$ for some $\hat c>0$.
Using $1 = \<n,n\> = \hat c^{\,2}-2\,\hat c\,\beta(N) +\|\beta\,\|^2_\alpha$, we get two values
\[
 \hat c = \beta(N)\pm(\beta(N)^2+c^2)^{1/2}.
\]
By condition $\beta(N)=0$ we have $\beta^\sharp\in W$, this yields $\hat c=c$ and \eqref{E-c-value}.
Thus,
\[
 \beta(n)=\beta(cN-\beta^\sharp)=-\|\beta\|^2_\alpha,\qquad
 1+\beta(n)=c^2.
\]
Finally, \eqref{E-c-value2} follows from \eqref{E-c-value0}.
\qed

\begin{lemma}\label{L-zZ}
 Let the Randers norm obeys $\,\beta(N)=0$ $($i.e., $\beta^\sharp\in W)$. If $u,U\in W$ and
\begin{equation}\label{E-u-Ucond}
 g(u,v)=\<U,v\>\quad  \mbox{\rm for all}\quad   v\in W
\end{equation}
then $\beta(u)=c^{-4}\beta(U)$ and
\begin{equation}\label{E-u-U}
 c^{2}\,u = U+c^{-2}\beta(U)\,\beta^\sharp.
\end{equation}
\end{lemma}

\proof
By \eqref{E-c-value2}, we have
\[
 g(u,v) = c^2 \<u -\beta(u) \beta^\sharp,\,v\>.
\]
Then from \eqref{E-u-Ucond}, since $u,U$ and $\beta^\sharp$ belong to $W$, we obtain
\[
 u -\beta(u)\beta^\sharp = c^{-2} U.
\]
Applying $\beta$ we get
 $\beta(u) -\beta(u)\,\|\beta\,\|^2_\alpha = c^{-2}\beta(U)$,
$\beta(u)=c^{-4}\beta(U)$ and then \eqref{E-u-U}.
\qed

\subsection{Finsler spaces}

Let $M^{m+1}$ be a connected smooth manifold and $TM$ its tangent bundle.
The natural projection $\pi:TM_0\to M$, where $TM_0:=TM\setminus\{0\}$ is called the \textit{slit tangent bundle}.
A~\textit{Finsler structure} on $M$ is a Minkowski norm $F$ in tangent spaces $T_pM$, which smoothly depends on a point $p\in M$.
 Note that $\pi_\ast$ maps the double tangent bundle $T^2M$ into $TM$ itself.

 A \textit{spray} on a manifold $M$ is a smooth vector field $\mathbb{G}$ on $TM_0$ such that
\begin{equation}\label{E-spray}
 \pi_\ast(\mathbb{G}_v) = v,\quad
 \mathbb{G}_{\lambda v} = \lambda \,(h_\lambda)_\ast(\mathbb{G}_v)\qquad (v\in TM_0,\ \lambda >0),
\end{equation}
where $h_\lambda : v\mapsto \lambda\,v$ is the homothety of $TM$.
 The meaning of \eqref{E-spray}$_1$ is that $\mathbb{G}$ is a \textit{second-order vector field} over $M$,
and \eqref{E-spray}$_2$ is the homogeneous quadratic condition.
In~local coordinates $(x^i)$, $\mathbb{G}$ is expressed as $\mathbb{G}(y)=y^i\partial_{x^i}-2 G^i\partial_{y^i}$,
where $G^i(\lambda\,y)=\lambda^2 G^i(y)$ $(\lambda>0)$.

Using $\mathbb{G}$ we define the follo\-wing notions: covariant derivative,
parallel translation (and parallel vectors) along a curve,
geodesics and curvature.
 A~curve $\gamma(t)$ in $TM_0$ satisfying $\dot\gamma=\mathbb{G}_\gamma$ is an integral curve of~$\mathbb{G}$;
it is equal to the canonical lift of $c:=\pi\circ\gamma$.
 The \textit{covariant derivative} of a vector field $u(t)$ along a curve $c(t)$ in $M$ is given by
 $D_{\dot c}\,u = \{\dot u^i + \Gamma^i_{kj}(\dot c)\,\dot c^k\,u^j\}\,\partial_{x^i\,|\,c}$\,.
Here $G^i=\frac12\,\Gamma^i_{kj}\,y^k y^j$ for smooth functions $\Gamma^i_{kj}=(G^i)_{y^k y^j}$ on $TM_0$,
see Corollary~\ref{L-hom-2}. The following properties are obvious:
\[
 D_{\dot c}\,(u+v)=D_{\dot c}\,u+D_{\dot c}\,v,\quad
 D_{\dot c}\,(f u)={\dot c}(f) \,u+f D_{\dot c}\,u,\quad
 D_{\lambda\dot c}\,u=\lambda D_{\dot c}\,u
\]
for any $f\in C^\infty(M)$ and $\lambda>0$, see \cite{sh2}.
A vector field $u(t)$ along $c$ is \textit{parallel} if $D_{\dot c}\,u(t)\equiv0$, i.e.,
\begin{equation*}
 \dot u^i + \Gamma^i_{kj}(\dot c)\,\dot c^k\,u^j=0\quad (i\ge1).
\end{equation*}
A curve $c(t)$ in $M$ is called a \textit{geodesic} of $\mathbb{G}$
if it is a projection of an integral curve of~$\mathbb{G}$; hence, $\ddot c=\mathbb{G}_{\dot c}$.
A curve $c(t)$ is a {geodesic} if and only if the tangent vector $u=\dot c$ is parallel along itself: $D_{\dot c}\,\dot c =0$.
For a geodesic $c(t)$ we have the following quasilinear system of second order~ODEs
\begin{equation*}
 \ddot c^{\,i} +2 G^{\,i}(\dot c) =0, \quad i = 1, \ldots, m+1\,.
\end{equation*}

A Finsler metric $F$ on $M$ induces a \textit{Finsler spray} ${\mathbb G}$ on $TM_0$,
whose geodesics are locally shortest paths connecting endpoints and have constant speed.
 Its {geodesic coefficients} are given by
\begin{equation*}
 G^{\,i} =\frac14\,g^{il}\big([F^2]_{x^k y^l}\,y^k -[F^2]_{x^l}\big)
 =\frac14\,g^{il}\big( 2\frac{\partial g_{jl}}{\partial x^k}
 -\frac{\partial g_{jk}}{\partial x^l}\big)\,y^j y^k\,,
\end{equation*}
see~\cite{sh2}.
Here $g_{ij}(y)=\frac12\,[F^2]_{y^i y^j}(y)$, compare \eqref{E-acsiom-M3}.
Then
 $\Gamma^i_{kj}(y)=\frac12\,g^{il}\big(\frac{\partial g_{jl}}{\partial x^k}
 +\frac{\partial g_{kl}}{\partial x^j} -\frac{\partial g_{jk}}{\partial x^l}\big)$
are homogeneous of 0-degree functions on $TM_0$.

\begin{remark}\rm
A Finsler metric on a manifold $M$ is called a \textit{Berwald metric} if in any local coordinate system $(x,y)$ in $TM_0$,
the Christoffel symbols $\Gamma^i_{jk}$ are functions on $M$ only,
in which case the geodesic coefficients $G^i=\frac12\,\Gamma^i_{kj}(x)\,y^k y^j$ are quadratic in $y=y^i\partial_{x^i}$.
On a {Berwald space}, the parallel translation along any geodesic preserves the Minkowski functionals;
thus, such spaces can be viewed as Finsler spaces modeled on a single Minkowski space.
Berwald metrics are characterized among Randers ones, $F = \alpha + \beta$,
by the following criterion: $\beta$ is parallel with respect to $\alpha$, see \cite[Theorem~2.4.1]{sh2}.
If~$\beta$ is a closed 1-form, then Finslerian geodesics are the same (as sets) as the geodesics of the metric $a$.
\end{remark}

A Finsler manifold is \textit{positively} (resp. \textit{negatively}) \textit{complete} if every geodesic $c(t)$ on $(0,t_0)$
can be extended for $(0,\infty)$ (resp. $(-\infty,0)$), and $F$ is
\textit{complete} if it is both positively and negatively complete.
This property is satisfied by all closed Finsler manifolds.
Let $(M,F)$ be positively complete;
 hence, for any $p,q\in M$ there exists a globally minimizing geodesic from $p$ to $q$,
see also Hopf-Rinov theorem \cite[p.~178]{sh2}.
Let $c_y$ be a geodesic with $c_y(0)=p$ and $\dot c_y(0)=y\in T_pM$.
The~\textit{exponential~map} is defined by $\exp_p(y)=c_y(1)$.
By~homogeneity of ${\mathbb G}$ one has $c_y(t)=c_{\,ty}(1)$ for $t>0$; hence, $\exp_p(ty)=c_y(t)$.
Recall \cite{sh1} that $\exp_p$ is smooth on $TM_0$ and $C^1$ at the origin with $d(\exp_p)_{|\,0}={\rm id}_{\,T_pM}$.

Consider a geodesic $c(t),\ 0\le t\le 1$. A $C^\infty$ map ${\cal H}:(-\eps,\eps)\times[0,1]\to M$ is called
a \textit{geodesic variation} of $c$ if ${\cal H}(0,t)=c(t)$ and for each $s\in(-\eps,\eps)$, the curve $c_s(t):={\cal H}(s,t)$ is a geodesic.
For a~geodesic variation ${\cal H}$ of $c$, the variation field $Y(t):=\frac{\partial{\cal H}}{\partial s}(0,t)$ along $c$
satisfies the \textit{Jacobi~equation}:
\begin{equation}\label{E-RC-Jacobi}
 D_{\dot c}D_{\dot c}\,Y +R_{\dot c}(Y) =0
\end{equation}
for some ($y\in TM$)-dependent (1,1)-tensor $R_y$. Jacobi equation
\eqref{E-RC-Jacobi} serves as the definition of curvature.
A vector field $Y(t)$ satis\-fying \eqref{E-RC-Jacobi} along a geodesic $c(t)$ is called \textit{Jacobi field}.
We have $g_{\dot c}(Y(t),\dot c(t))=\lambda^2(a+bt)$ and $g_{\dot c}(D_{\dot c}\,Y(t),\dot c(t))=\lambda^2b$
for some constants $a,b$ and $\lambda{=}F(\dot c)$.
The~orthogonal component $Y^\bot(t)=Y(t)-(a+bt)\dot c(t)$ of the Jacobi field $Y(t)$ along $c(t)$
is also a Jacobi field such that $Y^\bot(t)$ and $D_{\dot c}\,Y^\bot(t)$ are $g_{\dot c}$-orthogonal to $\dot c(t)$.
 Define $R^{(1)}_{\dot c(t)}:T_{c(t)}M\to T_{c(t)}M$ by $R^{(1)}_{\dot c(t)}\,(u(t))=D_{\dot c(t)}[R_{\dot c(t)}\,(u(t))]$,
where $u(t)$ is a parallel vector field along $c$.
Similarly, we define $R^{(2)}_{\dot c(t)}$, $R^{(3)}_{\dot c(t)}$ etc. Thus, by \eqref{E-RC-Jacobi}, a spray defines
transformations $R_y:T_pM\to T_pM$ called the \textit{Riemann curvature in a direction} $y\in T_pM\setminus\{0\}$, and
we have $R_y(y)=0$ and $R_{\,\lambda y}=\lambda^2 R_y\ (\lambda>0)$.
In~coordinates, $R_y=R^i_{\ k} dx^k\partial_{x_i}$ and $R^i_{\ k}(y)\,y^k=0$,
where $R^i_k$'s depend on the Finsler spray only~\cite{sh1}:
\begin{equation*}
 R^i_{\ k} = 2\,(G^i)_{x^k} -y^j\,(G^i)_{x^j\,y^k} +2\,G^j\,(G^i)_{y^j\,y^k} -(G^i)_{y^j}\,(G^j)_{y^k}\,.
\end{equation*}
Moreover,
$R^i_{\ k} =R^{\ i}_{j\ k l}\,y^j\,y^l$ for local
functions $\{R^{\ i}_{j\ k l}\}=\frac12\,(R^i_{\ k})_{y^j y^l}$ on $TM_0$ (see Corollary~\ref{L-hom-2}) and
\begin{eqnarray*}
 R^{\ i}_{j\ k l} = (\Gamma^i_{jl})_{x^k} -(\Gamma^i_{jk})_{x^l} +\Gamma^m_{jl}\,\Gamma^i_{mk} -\Gamma^m_{jk}\,\Gamma^i_{ml}\,.
\end{eqnarray*}
For the Finsler spray, $R_y$ is $g_y$-self-adjoint:
 $g_y(R_y(u),v)=g_y(u, R_y(v)),\ u,v\in T_pM$.

For a plane $P\subset T_pM$ tangent to $M$ and a vector $y\in P\setminus\{0\}$,
the \textit{flag curvature} $K(P,y)$ is~given~by
\[
 K(P,y)=\frac{g_y(R_y(u),u)}{g_y(y,y)g_y(u,u) - g_y(y,u)g_y(y,u)}\,,
\]
where $u\in P$ is such that $P={\rm span}\{y,u\}$; certainly, the value of $K(P,y)$  is independent of the choice of $u\in P$.
If $K(P,y)$ is a scalar function on $TM_0$ (that holds in dimension two) then $F$ is said to be of \textit{scalar} (flag) \textit{curvature},
in this case,
 $R_y(u)=K(\pi(y))\{g_y(y,y)u - g_y(y,y)y\}\ (y,u\in TM_0)$.
If~$K=K(\pi(y))$ (i.e., the flag curvature is \textit{isotropic}) and $m\ge2$ then $K=\const$,
see~\cite[Lemma 7.1.1]{cs2}.
For~each $K\in\RR$ there exist many non-isometric Finsler metrics of constant scalar curvature~$K$.

Let $\{e_i\}_{1\le i\le m+1}$ be a $g_y$-orthonormal basis for $T_pM$ such that $e_{m+1}=y/F(y)$,
and let $P_i={\rm span}\{e_i,y\}$ for some $y\in T_pM$.
Then $K(P_i,y)=F^{-2}(y)\,g_y(R_y(e_i),e_i)$.
The \textit{Ricci curvature} is a function on $TM_0$ defined as the trace of the Riemann curvature,
\[
 \Ric(y)=\sum\nolimits_{\,i=1}^m g_y(R_y(e_i),e_i)=F^2(y)\sum\nolimits_{i=1}^m K(P_i,y)
 \]
with the homogeneity property $\Ric(\lambda y)=\lambda^2\Ric(y)$\ $(\lambda>0)$.
In a coordinate system, by Corollary~\ref{L-hom-2} we have
$\Ric(y)=R^{\ i}_{j\ i k}\,y^j\,y^k=\Ric_{jk}\,y^j\,y^k$.
A Finsler space $(M^{m+1},F)$ is said to be of \textit{constant Ricci curvature} $\lambda$ (or, \textit{Einstein})~if
 $\Ric(y)=m\lambda\,F^2(y)\ (y\in TM_0)$,
or $\Ric_{jk}=m \lambda\,g_{jk}$ in coordinates.

\section{Codimension-one foliated Finsler spaces}
\label{sec:intform}

Given a transversally oriented codimension-one foliation $\calf$ of a Finsler manifold $(M^{m+1},F)$,
 there exists a globally defined $F$-normal (to the leaves) smooth vector field $n$
which defines a Riemannian metric $g:=g_n$ with the Levi-Civita connection~$\nabla$.
We have $g(n,u)=0\ (u\in T\calf)$ and $g(n,n)=F^2(n)$, see \eqref{E-c-value2c}.
 Then $\nu=n/F(n)$ is an $F$-unit normal.

\subsection{The Riemann curvature and the shape operator}

In this section we apply the variational approach to find a relationship between the Riemann curvature of $F$ and $g$. It generalizes the following.

\begin{proposition}[see \cite{sh2}]\label{P-R-Yt0}
Let $Y$ be a geodesic field on an open subset ${\cal U}$ in a Finsler space $(M,F)$
and $\hat g:=g_Y$ the induced metric on ${\cal U}$. Then the Riemann curvature of $F$ and $\hat F:=\sqrt{\hat g}$
obey $R_Y=\hat R_Y$. Moreover, $Y$ is a geodesic field of $\hat F$ and for the Levi-Civita connection we have $D_Y X=\hat D_Y X$.
\end{proposition}

For a codimension-one Riemannian foliation, a unit normal $\nu$ is a geodesic vector field;
hence, by Proposition~\ref{P-R-Yt0}, transformations $R_\nu$
defined for $F$ by \eqref{E-RC-Jacobi} coincide with the Jacobi operator $R(\cdot, \nu)\nu$ of the metric $g$.
Recall that the second differential is defined by $\nabla^2_{u,v}=\nabla_{u}\nabla_{v}-\nabla_{\nabla_u v}$ for any $u,v$.

\smallskip
Let $Y_t\ (|t|\le\eps)$ be a smooth family of $F$-unit vector fields on an open subset ${\cal U}$ in $(M,F)$.
Put~$\dot Y_t=\dt Y_t$ and $\dot g_t=\dt g_t$, where $g_t:=g_{\,Y_t}$ is a family of metrics on ${\cal U}$.
By definition  \eqref{E-Ctorsion} of the Cartan torsion, we have
\begin{equation}\label{E-gC}
 \dot g_t=2 C_{\,Y_t}(\,\cdot\,,\,\cdot\,,\dot Y_t).
\end{equation}
Note that $\dot g_t(Y_t,\,\cdot)=2 C_{\,Y_t}(\,Y_t,\,\cdot\,,\dot Y_t)=0$.

\begin{proposition}\label{P-R-Yt}
Let $Y_t\ (|t|\le\eps)$ doesn't depend on $t$ at a point $p\in{\cal U}$ and $u,v\in T_pM$. Then
\begin{eqnarray}\label{E-dt-Rm}
\nonumber
 -\dt R_t(u,Y_t,Y_t,v) \eq C_{\,Y}(u, \nabla^t_v Y_t, \nabla^t_{Y}\dot Y_t)
 +C_{\,Y}(\nabla^t_u Y_t, v,\nabla^t_{Y}\dot Y_t) \\
\nonumber
 \plus C_{\,Y}(\nabla^t_{Y}Y_t,v, \nabla^t_u\dot Y_t) +C_{\,Y}(u, \nabla^t_{Y}Y_t, \nabla^t_v\dot Y_t)\\
 \plus C_{\,Y}(u, v, (\nabla^t)^2_{Y,Y}\dot Y_t)
 +\,2(\nabla^t_{Y}C_{\,Y_t})(u, v, \nabla^t_{Y}\dot Y_t).
\end{eqnarray}
The shape operators $A_t$ (when $Y_p=\nu_p$) of $\calf$ with respect to $g_t$ and the volume forms ${\rm d}V_t$
at $p$ obey
\begin{eqnarray}\label{E-dt-A}
 g_t(\dt A_t(u),v) \eq -C_{\,\nu}(u, v, \nabla^t_{\nu}\,\dot Y_t),\quad
 \dt ({\rm d}V_t) = 0.
\end{eqnarray}
\end{proposition}

\proof
Put $\Pi(u,v)=\dt\nabla^t_{u}\,v$ for $t$-independent vector fields $u,v$.
Then, see \cite{topp},
\begin{equation}\label{E-Pi}
 2\,g_t(\Pi(u,v),w) = (\nabla^t_v\,\dot g_t)(u,w) +(\nabla^t_u\,\dot g_t)(v,w) -(\nabla^t_w\,\dot g_t)(u,v),
\end{equation}
and for arbitrary $t$-dependent vector fields $X_t$ and $Z_t$ we obtain
\[
 \dt\nabla^t_{X_t} Z_t = \Pi(X_t,Z_t) +\nabla^t_{X_t}(\dt Z_t) +\nabla^t_{\dt X_t} Z_t.
\]
 By definition,
\[
 R_t(u,Z_t)Y_t=\nabla^t_u(\nabla^t_{Z_t} Y_t)-\nabla^t_{Z_t}(\nabla^t_u Y_t)-\nabla^t_{[u, Z_t]}Y_t.
\]
So,
\[
 \dt R_t(u,Z_t)Y_t=\dt(\nabla^t_u(\nabla^t_{Z_t}\,Y_t))-\dt(\nabla^t_{Z_t}(\nabla^t_u\,Y_t))-\dt(\nabla^t_{[u, Z_t]}\,Y_t).
\]
Deriving the terms of the above,
\begin{eqnarray*}
 \dt(\nabla^t_{Z_t}(\nabla^t_u\,Y_t)) \eq \Pi(Z_t,\nabla^t_u\,Y_t)
 +\nabla^t_{Z_t}(\Pi(u,Y_t)) +\nabla^t_{Z_t}(\nabla^t_u\,\dot Y_t) +\nabla^t_{\dot Z_t}(\nabla^t_u\,Y_t),\\
 \dt(\nabla^t_u(\nabla^t_{Z_t}\,Y_t)) \eq
 \Pi(u,\nabla^t_{Z_t}\,Y_t) +\nabla^t_u(\Pi(Z_t, Y_t)) +\nabla^t_u(\nabla^t_{\dot Z_t}\,Y_t)
 +\nabla^t_u(\nabla^t_{Z_t}\,\dot Y_t),\\
 \dt(\nabla^t_{[u, Z_t]}\,Y_t) \eq \Pi([u,Z_t],Y_t)+\nabla^t_{[u,Z_t]}\,\dot Y_t+\nabla^t_{[u,\dot Z_t]}\,Y_t
\end{eqnarray*}
with $\dot Z_t=\dt Z_t$, we obtain a `time-dependent' version of \cite[Proposition~2.3.4]{topp},
\begin{eqnarray*}
 \dt R_t(u,Z_t)Y_t = (\nabla^t_u\,\Pi)(Z_t,Y_t) -(\nabla^t_{Z_t}\,\Pi)(u,Y_t)
 +R_t(u,Z_t)\dot Y_t +R_t(u,\dot Z_t)Y_t.
\end{eqnarray*}
We shall compute $\dt R_t(u, Y_t,Y_t,v):=\dt g_t(R_t(u, Y_t)Y_t,v)$ at $p$\,; thus, terms with $\dot Y$ will be canceled at the final~stage.
Assume at a `time' $t$ of our choice, $\nabla=\nabla^t$ and $\nabla u=\nabla v=0$ at $p$.
Then perform the following preparatory calculations at~$p$\,:
\begin{eqnarray*}
 \frac12\,Y\big((\nabla^t_u\, \dot g_t)(Y_t,v)\big)
 \eq Y\big(u\,(C_{\,Y_t}(Y_t,v,\dot Y_t)) -C_{\,Y_t}(\nabla^t_u\,Y_t,v,\dot Y_t)\big)\\
  \eq -C_{\,Y}(\nabla_u Y_t, v, \nabla_{Y}\dot Y_t) ,\\
 \frac12\,Y\big((\nabla^t_{Y_t}\, \dot g_t)(u,v)\big)
  \eq Y\big({Y_t}\,(C_{\,Y_t}(u,v,\dot Y_t))\big)
  -Y(C_{\,Y_t}(\nabla^t_{Y_t}\,u,v,\dot Y_t)) \\
  && -\,Y(C_{\,Y_t}(u,\nabla^t_{Y_t}\,v,\dot Y_t)) \\
 \eq C_{\,Y}(u, v, \nabla_{Y}\nabla_{Y_t}\,\dot Y_t) + 2(\nabla_{Y} C_{\,Y})(u, v, \nabla_{Y}\dot Y_t),\\
 \frac12\,Y\big((\nabla^t_v\, \dot g_t)(u,Y_t)\big)
 \eq Y\big(v\,(C_{\,Y_t}(u,Y_t,\dot Y_t)) -C_{\,Y_t}(u,\nabla_v Y_t,\dot Y_t)\big)\\
 \eq -C_{\,Y}(u,\nabla_v\,Y_t,\nabla_{Y}\,\dot Y_t),\\
 (\nabla_{\nabla_{Y}Y_t}\,\dot g_t)(u,v) \eq 2C_{\,Y}(u, v, \nabla_{\nabla_{Y}Y_t}\dot Y_t),\\
 (\nabla_u\,\dot g_t)(\nabla_{Y}Y_t,v) \eq 2C_{\,Y}(\nabla_{Y}Y_t, v, \nabla_u\dot Y_t), \\
 (\nabla_v\,\dot g_t)(u,\nabla_{Y}Y_t) \eq 2C_{\,Y}(u, \nabla_{Y}Y_t, \nabla_v\dot Y_t)\,.
\end{eqnarray*}
Using all of that and \eqref{E-gC} we obtain at $p$:
\begin{eqnarray*}
 &&\hskip-10mm \<(\nabla_{Y}\,\Pi)(u,Y_t), v\> = \<\nabla_{Y}\,(\Pi(u,Y_t)) - \Pi(u,\nabla_{Y}\,Y_t),\, v\>\\
 \eq Y\<\Pi(u,Y_t),\,v\> - \<\Pi(u,\nabla_{Y}\,Y_t),\, v\>\\
 \eq \frac12\,Y\big[\,(\nabla^t_u\,\dot g_t)(Y_t,v) +(\nabla^t_{Y_t}\,\dot g_t)(u,v) -(\nabla^t_v\,\dot g_t)(u,Y_t)\,\big]\\
 && -\,\frac12\,\big[\,
 (\nabla_{\nabla_{Y}Y_t}\,\dot g_t)(u,v) +(\nabla_u\,\dot g_t)(\nabla_{Y}Y_t,v) -(\nabla_v\,\dot g_t)(u,\nabla_{Y}Y_t)
 \,\big]\\
 \eq C_{\,Y}(u, \nabla_v Y_t, \nabla_{Y}\dot Y_t) -C_{\,Y}(\nabla_u Y_t, v,\nabla_{Y}\dot Y_t) \\
 && +\,2(\nabla_{Y}\,C_{\,Y_t})(u, v, \nabla_{Y}\dot Y_t) +C_{\,Y}(u, v, \nabla_{\,Y}\nabla^t_{\,Y_t}\dot Y_t) -C_{\,Y}(u, v, \nabla_{\nabla_{Y}Y_t}\dot Y_t) \\
 && -\,C_{\,Y}(\nabla_{Y}Y_t, v, \nabla_u\dot Y_t) +C_{\,Y}(u, \nabla_{Y}Y_t, \nabla_v\dot Y_t).
\end{eqnarray*}
Here the terms with $C_{\,Y}(Y,\,\cdot\,,\cdot\,)$ were canceled on ${\cal U}$, and
the identity $[Y_t,v]^\top = -(\nabla^t_v\,Y_t)^\top$
at $p$ (where $^\top$ is the orthogonal to $Y$ at $p$ component of a vector) was applied.
Similarly, we use at $p$
\begin{eqnarray*}
 && u\big[(\nabla^t_{Y_t}\,\dot g_t)(Y_t,v)\big]= -2C_{\,Y}(\nabla_{Y}Y_t, v, \nabla_u\dot Y_t),\quad
 u\big[(\nabla^t_v\,\dot g_t)(Y_t,Y_t)\big] = 0,\\
 && (\nabla_{\nabla_u Y_t}\,\dot g)(Y,v)=0,\quad
 (\nabla_v\,\dot g)(Y,\nabla_u Y_t)=0,\\
 && (\nabla_{Y}\,\dot g)(\nabla_u Y_t,v)=2\,C_{\,Y}(\nabla_u\,Y_t, v, \nabla_{Y}\dot Y_t)
\end{eqnarray*}
to find
\begin{eqnarray*}
 &&\<(\nabla_u\,\Pi)(Y_t,Y_t), v\> = \<\nabla_{u}(\Pi(Y_t,Y_t)) - 2\Pi(Y_t,\nabla_u Y_t), v\>\\
 && = u\<\Pi(Y_t,Y_t), v\> - 2\,\<\Pi(Y_t,\nabla_u Y_t), v\>\\
 && = u\big[(\nabla^t_{Y_t}\,\dot g_t)(Y_t,v) -\frac12\,(\nabla^t_v\,\dot g_t)(Y_t,Y_t)\big]\\
 &&\quad -(\nabla_{\nabla_u Y_t}\,\dot g)(Y_t,v)-(\nabla_{Y}\,\dot g)(\nabla_u Y_t,v) +(\nabla_v\,\dot g)(Y,\nabla_u Y_t)\\
 && =-2C_{\,Y}(\nabla_{Y}Y_t,v, \nabla_u\dot Y_t) -2\,C_{\,Y}(\nabla_u Y_t, v, \nabla_{Y}\dot Y_t).
\end{eqnarray*}
Since $\dot Y=0$ at $p$, we have
\begin{eqnarray*}
 && \dt R_t(u,Y_t,Y_t,v) = (\dt g)(R_t(u,Y_t)Y_t,v)+ g(\dt R_t(u,Y_t)Y_t,v)\\
 && = 2\,C_Y(R_t(u,Y_t)Y_t,v, \dot Y)+ g(\dt R_t(u,Y_t)Y_t,v)=g(\dt R_t(u,Y_t)Y_t,v).
\end{eqnarray*}
Finally, we have \eqref{E-dt-Rm} at $p$ for all $t\ge0$.
For the second fundamental form $b_t$ of $\calf$ (with respect to~$g_t$),
as in the proof of \cite[Lemma~2.9]{rw1}, using \eqref{E-gC}, \eqref{E-Pi}, $\dot g(p)=0$ and $\dot Y(p)=0$,
we~get at a point $p$:
\begin{eqnarray*}
 \dt b_t(u,v) \eq \dot g(\nabla_u v, Y)+g(\dt\nabla_u v, Y)+g(\nabla_u v, \dt Y)\\
 \eq \frac12\,\big((\nabla_u\dot g)(v,Y) +(\nabla_v\dot g)(u,Y) -(\nabla_Y\dot g)(u,v)\big) +g(\nabla_u v,\dot Y)  \\
 \eq -\nabla_{Y} (C_Y(u, v, \dot Y)) = -C_Y(u, v, \nabla_{Y}\dot Y).
\end{eqnarray*}
 From this, using $b_t(u,v)=g_t(A_t(u),v)$, we get \eqref{E-dt-A}$_1$:
\begin{eqnarray*}
 g_t(A_t(u),v) = \dt b_t(u,v) -\dot g(A(u),v) =-C_\nu(u, v, \nabla_{\nu}\,\dot Y).
\end{eqnarray*}
By the formula for the volume form of a $t$-dependent metric, $\dt ({\rm d}V_t) =\frac12\,(\tr\dot g)\,{\rm d}V_t$, see \cite{topp}, and definition of the mean Cartan torsion, we~get
\begin{equation}\label{E-dtV}
 \dt ({\rm d}V_t) = I_{Y_t}(\dot Y_t)\,{\rm d}V_t.
\end{equation}
Next, \eqref{E-dt-A}$_2$ follows from \eqref{E-dtV} and $\dot Y(p)=0$.
\qed

\smallskip
Let $L$ be a leaf through a point $p\in M$, and $\rho$ the local distance function to $L$ in a neighborhood of~$p$.
Denote by $\hat\nabla$ the Levi-Civita connection of the (local again) Riemannian metric $\hat g:= g_{\,\nabla\rho}$.
Note that $\nabla\rho=\nu$ on~$L$.
 The~\textit{shape operator} $A:T\calf\to T\calf$ (self-adjoint for $g$) is defined at $p\in M$ by (compare \cite{sh2} with the opposite sign)
\[
 A(u)=-{\hat\nabla_u\,\nu}\quad (u\in T_p\calf).
\]
The shape operator $A^g:T\calf\to T\calf$ with respect to the metric $g$ is defined at $p\in M$~by
\[
 A^g(u)=-{\nabla_u\,\nu}\quad (u\in T_p\calf).
\]
Note that $2\,g(\nabla_u\,\nu,\nu)=u(g(\nu,\nu))=0\ (u\in T\calf)$; hence, $\nabla_u\,\nu\in T\calf$.
The \textit{mean curvature} function (of the leaves with respect to $g$) is defined by $H^g=\tr A^g$.
Recall that $\calf$ is $g$-\textit{totally umbilical} if $A^g=H^g I_m$, and
is $g$-\textit{totally geodesic} if $A^g\equiv 0$.

\begin{corollary}\label{C-R-Yt}
Let $L$ be a hypersurface in an open set ${\cal U}\subset M$.
If an $F$-unit vector field $Y_t\ (0\le t\le\eps)$ is given in ${\cal U}$ and orthogonal to $L$
then for the metric $g_t:=g_{\,Y_t}$ for all $u,v\in T_pL\ (p\in L)$ we~have
\begin{eqnarray}\label{E-dt-Rm-L}
\nonumber
 \dt R_t(u,Y_t,Y_t,v) \eq C_{\,Y}(A_t(u), v,\nabla^t_{Y}\dot Y_t) +C_{\,Y}(u, A_t(v),\nabla^t_{Y}\dot Y_t)\\
 && -\,C_{\,Y}(u, v, (\nabla^t)^2_{\,Y,Y}\dot Y_t) - 2(\nabla^t_{Y} C_{\,Y_t})(u, v, \nabla^t_{Y}\dot Y_t),\\
\label{E-dt-Ab}
 g(\dt A_t(u),v) \eq -C_{\,Y}(u, v, \nabla^t_{Y}\dot Y_t),\quad
 \dt ({\rm d}V_t) = 0.
\end{eqnarray}
\end{corollary}

\proof This follows from $\dot Y_t=0$ on $L$, the definition of
$A_t$ (for $g_t$) and \eqref{E-dt-Rm}--\eqref{E-dt-A}.
\qed

\begin{definition}\rm
 A vector field $\widehat Y$ defined in some neighborhood ${\cal U}\subset M$ of a point $p\in{\cal U}$ is called
 a~\textit{geodesic extension} of a vector $Y_p\in T_pM$ if $\widehat Y(p)=Y_p$
 and the integral curves of $\widehat Y$ are geodesics of the Finsler metric.
 Similarly, we define a~\textit{geodesic extension} of a (e.g. normal) vector field along a hypersurface $L\subset{\cal U}$.
 In both cases, $\hat g:=g_{\widehat Y}$ is called the \textit{osculating Riemannian metric} of $F$ on ${\cal U}$.
\end{definition}

We will use osculating metric (given locally) to express the Riemannian curvature of $g=g_\nu$ (for an unit $F$-normal $\nu$ to $\calf$) in terms of Riemannian curvature and the Cartan torsion of $F$.

Given a vector field $Y\in C^\infty(TM)$, let $C_Y^\sharp$ be a $(1,1)$-tensor $g_{\,Y}$-dual to the symmetric bilinear form $C_{\,Y}(\cdot\,,\cdot\,,\nabla_Y\,Y)$.
Note that
$C_n(\cdot\,,\cdot\,,\nabla_n\,n)=C_{c^2\nu}(\cdot\,,\cdot\,,c^4\nabla_\nu\,\nu)=c^2C_{\nu}(\cdot\,,\cdot\,,\nabla_\nu\,\nu)$.

\begin{theorem}\label{T-unit-n}
Let $\nu$ be a unit normal to a codimension-one foliation of a Finsler space $(M^{m+1},F)$.
The Riemann curvatures $($in the $\nu$-direction$)$ of $F$ and $g=g_\nu$ are related by
\begin{eqnarray}\label{E-dt-Rm-Rm}
\nonumber
 && g((R_{\,\nu} - R^g_{\,\nu})(u),v) = -C_{\nu}\big(A^g(u)+\frac12\,C_\nu^\sharp(u), v, \nabla_{\nu}\,\nu\big)\\
\nonumber
 &&\hskip10mm -\,C_{\nu}\big(u, A^g(v)+\frac12\,C_\nu^\sharp(v), \nabla_{\nu}\,\nu\big)\\
 &&\hskip10mm +\,C_{\nu}\big(u, v, \nabla^2_{\nu,\nu}\,\nu -C_\nu^\sharp(\nabla_{\nu}\,\nu)\big)
 +\,2(\nabla_{\nu}C_{\nu})(u,v,\nabla_{\nu}\,\nu) \quad (u,v\in T_pL).
\end{eqnarray}
The shape operators and volume forms are related by
\begin{equation}\label{E-dt-A2}
  A - A^g = C_\nu^\sharp,\qquad
  {\rm d}V_g = e^{\tau(\nu)}\,{\rm d}V_F.
\end{equation}
In particular, the traces are related by
\begin{eqnarray}\label{E-dt-Ric-Ric}
\nonumber
 \Ric_{\nu}-\Ric^g_{\,\nu} \eq
 I_{\nu}(\nabla^2_{\nu,\nu}\,\nu -C_\nu^\sharp(\nabla_{\nu}\,\nu)) +2(\nabla_{\nu}\,I_{\nu})(\nabla_{\nu}\,\nu)\\
 && -\,\tr\big(C^\sharp_{\,\nu}( C_\nu^\sharp + 2\,A^g)\big),\\
 \nonumber
 \tr A-\tr A^g \eq I_{\nu}(\nabla_{\nu}\,\nu).
\end{eqnarray}
\end{theorem}

\proof
Let ${\cal U}$ be a ``small" neighborhood of $p\in L$
such that any two geodesics starting from $L\cap{\cal U}$ in the $\nu$-direction do not intersect in~${\cal U}$.
Then for any $q\in{\cal U}$ there is a unique
geodesic $\gamma$ starting from $L$ in the $\nu$-direction
such that $\gamma(s)=q$ for some $s\ge0$, in other words, $q=\exp_{\gamma(0)}(s\,\dot\gamma(0))$.
Thus, $\widehat Y: q \to\dot\gamma(s)\ (q\in{\cal U})$ is an $F$-unit geodesic vector field
($\nabla_{\widehat Y}\widehat Y=0$) -- a geodesic extension of $\nu_{\,|\,L}$.

Consider a family of vector fields $Y_t=t\,\widehat Y +(1-t)\,\nu\ (0\le t\le 1)$ on~${\cal U}$, define the Riemannian metrics $g_t:=g_{\,Y_t}$,
$g_1$ being osculating, and denote by $R_t$ their Riemann curvatures.
Since $\dot Y_t=\widehat Y-\nu$ and $Y_{t\,|L}=\nu_{\,|L}=\widehat Y_{\,|L}$ for all~$t$,
we have $\dot Y_{t\,|L}=0$ and $g_{t\,|L}\equiv g_{\,|L}$.
 By \eqref{E-gC} and \eqref{E-Pi}, we get $\Pi_t(\nu,\nu)=\Pi_t(\nu,\widehat Y)=0$ on $L$; hence,
$\nabla^t_\nu\,\nu$ and $\nabla^t_{\nu}\widehat Y$ restricted on $L$ don't depend on $t$.
Next,~we~find
\[
 g(\Pi(\nu,\nu),v)=C_\nu(u,v,\nabla_\nu(\widehat Y-\nu))=-C_\nu(u,v,\nabla_\nu\,\nu),\quad u,v\in TM_{\,|L},
\]
i.e., $\Pi(\nu,u)=-C^\sharp_\nu(u)$. We calculate on $L$:
\begin{eqnarray*}
 g(\dt(\nabla^t_\nu\,u), v) \eq \nabla^t_\nu(C_Y(u,v,\widehat Y-\nu)) +\nabla^t_u(C_Y(\nu,v,\widehat Y-\nu))
 -\nabla^t_v(C_Y(u,\nu,\widehat Y-\nu))\\
 \eq (\nabla^t_\nu C_Y)(u,v,\widehat Y-\nu) +C_Y(u,v,\nabla^t_\nu(\widehat Y-\nu))\\
 \plus (\nabla^t_u C_\nu)(n,v,\widehat Y-\nu) +C_\nu(\nabla^t_u\,\nu,v,\widehat Y-\nu)
 +C_\nu(\nu,v,\nabla^t_u(\widehat Y-\nu))  \\
 \minus (\nabla^t_v C_\nu)(u,\nu,\widehat Y-\nu) - C_\nu(u,\nabla^t_v\,\nu,\widehat Y-\nu) -C_\nu(u,\nu,\nabla^t_v(\widehat Y-\nu)) \\
 \eq C_\nu(u,v,\nabla^t_\nu(\widehat Y-\nu)) = -C_\nu(u,v,\nabla_\nu\,\nu).
\end{eqnarray*}
Since,
$\dt(g(\nabla^t_\nu\,u,v))=g(\dt\nabla^t_\nu\,u,v)$ and $\dt(g(\nabla^t_u\,\nu,v))=g(\dt\nabla^t_u\,\nu,v)$ on $L$,
we obtain
\begin{eqnarray*}
 g(\nabla^t_\nu\,u, v)=g(\nabla_\nu\,u,v) -t\,C_\nu(u,v,\nabla_\nu\,\nu),\\
 g(\nabla^t_u\,\nu, v)=g(\nabla_u\,\nu,v) -t\,C_\nu(u,v,\nabla_\nu\,\nu).
\end{eqnarray*}
Recall that $\nabla^2_{u,v}$ is tensorial in $u,v$.
We show that $(\nabla^t)^2_{\nu,\nu}\,\widehat Y$ is $t$-independent on $L$:
\begin{eqnarray*}
 (\nabla^t)^2_{\widehat Y,\widehat Y}\,\widehat Y \eq \nabla^t_n(\nabla^t_{\widehat Y}\,\widehat Y)
 = \nabla_\nu(\nabla^t_{\widehat Y}\,\widehat Y) -t\,C^\sharp_\nu(\nabla^t_{\nu}\,\widehat Y)\\
 \eq \nabla_\nu(\nabla^t_{\widehat Y}\,\widehat Y)
 =\nabla_\nu(\nabla_{\widehat Y}\,\widehat Y -t\,C^\sharp_\nu(\widehat Y))\\
 \eq \nabla^2_{\nu,\nu}\,\widehat Y
 -t\,(\nabla_\nu C^\sharp_\nu)(\widehat Y) - t\,C^\sharp_\nu(\nabla_\nu\widehat Y)
 =\nabla^2_{\nu,\nu}\,\widehat Y.
\end{eqnarray*}
Thus, $(\nabla^2_{\nu,\nu}\,\widehat Y)_{\,|L}=(\widehat\nabla^2_{\nu,\nu}\,\widehat Y)_{\,|L}=0$.
Using this and $(\nabla_\nu\widehat Y)_{\,|L}=0$, we find on~$L$:
\begin{eqnarray*}
 \nabla^t_{Y_t}\,\dot Y_t \eq -\nabla_\nu\,\nu,\\
 (\nabla^t)^2_{Y_t,Y_t}\,\dot Y_t \eq (\nabla^t)^2_{\nu,\nu}\,(\widehat Y-\nu)
 =\nabla^t_\nu\big(\nabla_\nu(\widehat Y-\nu) -t\,C^\sharp_\nu\,(\widehat Y-\nu)\big) \\
 \eq \nabla^2_{\nu,\nu}\,(\widehat Y-\nu) -t\,\nabla_\nu\,(C^\sharp_\nu(\widehat Y-\nu)) -t\,C^\sharp_\nu(\nabla_\nu\,(\widehat Y-\nu))\\
 \eq -\nabla^2_{\nu,\nu}\,\nu +2\,t\,C_\nu^\sharp(\nabla_\nu\,\nu).
\end{eqnarray*}
Then we obtain on $L$:
\begin{eqnarray*}
 C_{\,Y_t}(\cdot\,,\cdot\,,\nabla_{Y_t}\dot{\,Y_t}) \eq C_{\,\nu}(\cdot\,,\cdot\,,\nabla_{\nu}(\widehat Y-\nu))
 = -C_{\,\nu}(\cdot\,,\cdot\,,\nabla_{\nu}\,\nu),\\
 C_{\,Y_t}(\cdot\,,\cdot\,,\nabla^2_{\,Y_t,Y_t}\dot{\,Y_t}) \eq C_{\,\nu}(\cdot\,,\cdot\,,\nabla^2_{\,\nu,\nu}(\widehat Y-\nu))
 = -C_{\,\nu}(\cdot\,,\cdot\,,\nabla^2_{\,\nu,\nu}\,\nu).
\end{eqnarray*}
Next, we calculate on $L$, using $C_{Z}(Z,\cdot\,,\cdot\,)=0$ for $Z=\nabla_\nu\,\nu$,
\begin{eqnarray*}
 && (\nabla_{Y_t} C_{\,Y_t})(\cdot\,,\cdot\,,\nabla_{Y_t}\dot{\,Y_t})
 =(\nabla_\nu C_{\,t\,\widehat Y+(1-t)\,\nu})(\cdot\,,\cdot\,,-\nabla_\nu\,\nu) \\
 && =(\nabla_\nu C)_{\,\nu}(\cdot\,,\cdot\,,-\nabla_\nu\,\nu) +C_{\,(1-t)\nabla_\nu\,\nu}(\cdot\,,\cdot\,,-\nabla_\nu\,\nu)
 =-(\nabla_\nu C_{\,\nu})(\cdot\,,\cdot\,,\nabla_\nu\,\nu).
\end{eqnarray*}
By the above and \eqref{E-dt-A}$_1$, we obtain \eqref{E-dt-A2}$_1$.
By Corollary~\ref{C-R-Yt}, for all $t\in[0,1]$, and using
$A_t=A^g+t\,C^\sharp_\nu$, see \eqref{E-dt-A2}$_1$, and $(\nabla^t)^2_{\nu,\nu}\,\nu=-\nabla^2_{\nu,\nu}\,\nu +2\,t\,C_\nu^\sharp(\nabla_{\nu}\,\nu)$, we obtain
\begin{eqnarray*}
 \dt R_t(u,\nu,\nu,v) \eq - C_{\nu}(A_t(u), v, \nabla_{\nu}\,\nu) -C_{\nu}(u, A_t(v), \nabla_{\nu}\,\nu)\\
 && +\,C_{\nu}(u, v, (\nabla^t)^2_{\nu,\nu}\,\nu) +2(\nabla_{\nu}C_{\nu})(u, v, \nabla_{\nu}\,\nu)\\
 \eq -C_{\nu}(A^g(u)+t\,C^\sharp_\nu(u), v, \nabla_{\nu}\,\nu) - C_{\nu}(u, A^g(u)+t\,C^\sharp_\nu(v), \nabla_{\nu}\,\nu)\\
 && +\,C_{\nu}(u, v, -\nabla^2_{\nu,\nu}\,\nu +2\,t\,C_\nu^\sharp(\nabla_{\nu}\,\nu))
 +\,2(\nabla_{\nu}C_{\nu})(u, v, \nabla_{\nu}\,\nu)
\end{eqnarray*}
for $u,v\in T_pL$, where the right hand side becomes linear in $t$.
Integrating this by $t\in[0,1]$ yields \eqref{E-dt-Rm-Rm}.
Finally, using the equality for volume forms, ${\rm d}\widehat V = {\rm d}V_g$, and definition of $\tau$ (see Section~\ref{sec:M-R-norms}),
we get \eqref{E-dt-A2}$_2$.
\qed

\smallskip
Since any geodesic vector field $Y$ satisfies conditions
\begin{equation}\label{E-CYY}
  C_Y(u,v,\nabla_Y\,Y) = 0,\quad C_Y(u,v,\nabla^2_{Y,Y}\,Y) = 0\quad (\forall\,u,v),
\end{equation}
the following corollary generalizes Proposition~\ref{P-R-Yt0}.

\begin{corollary}\label{C-unit-n}
If $Y$ is a unit vector field on a Finsler space $(M,F)$ and $g:=g_Y$ a Riemannian metric on $M$
with the Levi-Civita connection $\nabla$ and conditions \eqref{E-CYY}, then $R_Y=R^g_{\,Y}$.
\end{corollary}

\proof
By \eqref{E-CYY}, we have $C^\sharp_Y=0$ and
\[
 (\nabla_Y C_Y)(u,v,\nabla_Y\,Y)=\nabla_Y (C_Y(u,v,\nabla_Y\,Y))-C_Y(u,v,\nabla^2_{Y,Y}\,Y)=0.
\]
 If a vector field $\,\widehat Y$
is a local geodesic extension of $Y(p)$ then $R^g_{\,Y}=\hat R_Y$ (and $A^g = \hat A$) at $p$, see \eqref{E-dt-Rm-Rm} and \eqref{E-dt-A2}.
 Thus, the claim follows from Proposition~\ref{P-R-Yt0}.
\qed

\subsection{Integral formulae}
\label{sec:IF}

Let $\cal F$ is a codimension-one foliation of a closed
 Finsler space $(M^{m+1}, F)$ with the Busemann-Hausdorff volume form ${\rm d}V_F\,$.
Define a family of diffeomorphisms $\{\phi_t: M\to M,\ 0\le t<\eps\}$ ($\varepsilon>0$ being small enough) by
\begin{equation*}
 \phi_t(p)=\exp_p(t\,\nu),\quad
 {\rm where}\quad \nu\in T_pM\quad\mbox{is an $F$-unit normal to}\ \calf \ {\rm at}\ p\in M.
\end{equation*}
Let $c(t)\ (t\ge0)$ be an $F$-geodesic with $c(0)=p$ and $\dot c(0)=\nu(p)$.
Any geodesic variati\-on built of $\phi_t$-trajectories determines an $F$-Jacobi field $Y(t)$ on $c$, and
$A_p(Y(0))=-[D_{\dot c(t)}\,Y(t)]_{\,|t=0}$, see~\cite[p.~225]{sh2}.
Recall that if vectors $u(t)$ and $v(t)$ are $D$-parallel along
$c(t)$ then $g_{\dot c(t)}(u(t), v(t))$ is constant.
Choose a positively oriented $g_{\nu(p)}$-ortho\-normal frame $(e^1,\dots, e^m)$ of
$T_p\calf$ and extend it by parallel translation to the frame
$(E_t^1,\dots, E_t^m)$ of vector fields $g_{\dot c(t)}$-orthogonal to $\dot c(t)$ along $c(t)$.
Denote also by $E_t^{m+1}=\dot c(t)$ the tangent vector field along $c(t)$.
Denote by $Y^i(t)\ (i\le m)$ the Jacobi field along $c(t)$ satisfying
$Y^i(0)=e^i$ and $D_{\dot c}\,Y^{i}(0)=A_p(e^i)$.
Let $R(t)$ be the matrix with entries $g_{\dot c}(R_{\dot c}(E_t^i), E_t^j)$.
Denote by ${\bf Y}(t)$ the $m\times m$ matrix consisting of the
scalar products $g_{\dot c}(Y^i(t), E_t^j)$ (``$F$-{\em Jacobi tensor}").
Then~${\bf Y}(0)=I_{m}$ and ${\bf Y}^{\prime}(0)=A_p$.
It is known (see, for instance, \cite[Sections 2.1 and 2.2]{sh2}) that
\begin{equation*}
 |\,d\,\phi_t(p)|=\det{\bf Y}(t)\,,
\end{equation*}
where $|\,d\phi_t(p)|$ is the Jacobian of $\phi_t$ at $p$.
Assume that $R^{(1)}_{\dot c(t)}\equiv0$ for any $F$-geodesic $c(t)\ (t\ge0)$
(e.g. $(M, F)$ is \textit{locally symmetric} with respect to $F$).
For $t=0$, we have $R^{(2)}_{\dot c(0)}\equiv R^{(3)}_{\dot c(t)}\equiv\ldots\equiv0$.
 For short, write $ R_p:=R(0)$.
 Note that ${\tr}\,R_p=\Ric(\nu(p))$.
The $F$-Jacobi equation ${\bf Y}^{\prime\prime}=-R(t){\bf Y}$ implies~that
\begin{equation*}
 {\bf Y}^{(2k  )}(0)=(-R_p)^k,\quad
 {\bf Y}^{(2k+1)}(0)=(-R_p)^k A_p,
 \quad k=0,1,2,\dots
\end{equation*}
Hence, our Jacobi tensor has the form
\[
 {\bf Y}(t)=\sum\nolimits_{\,k=0}^\infty {\bf Y}^{(k)}(0)\,\frac{t^k}{k!}
 =I_{m}+t A_p-\frac{t^2}{2!} R_p-\frac{t^3}{3!} R_p A_p+\frac{t^4}{4!} R_p^2
 +\ldots.
\]
Certainly, the radius of convergence of the series is uniformly bounded from below on $M$ (by $1/\|R\|_F>0$).
The volume of $M$ is defined by ${\rm Vol}_F(M)=\int_M\,{\rm d}V_F\,$.
Therefore -- by Domina\-ted Convergence Theorem -- its integration
together with Change of Variables Theorem yield the equa\-lity for any $t\ge0$ small enough
\begin{equation}\label{eq5x}
 {\rm Vol}_F(M)=\int_M\det\big(I_{m}+t A_p-\frac{t^2}{2!} R_p
 -\frac{t^3}{3!}  R_p A_p + \frac{t^4}{4!}  R_p^2 +\ldots\big)\,{\rm d}V_F\,,
\end{equation}
where ${\rm d}V_F$ is the volume form of $F$.
Formula (\ref{eq5x}) together with Lemma~\ref{lem2} of Appendix imply our main result
(which generalizes that of \cite{rw2} valid for the Riemannian case).
Note that the invariants $\sigma_{{\lambda}} (A_1, \ldots , A_k)$ of a set of real $m\times m$ matrices $A_i$ are defined and discussed in Appendix.

\begin{theorem}\label{thm:main}
 If $\cal F$ is a codimension-one foliation on a closed Finsler manifold $(M^{m+1}, F)$,
 which is $F$-locally symmetric, then for any $0\le k\le m$ one has
\begin{equation}\label{eq:main}
 \int_M \sum\nolimits_{\,\|{\lambda}\|=k}\sigma_{{\lambda}}\left(B_1(p), \ldots B_k(p)\right)\,{\rm d}V_F=0,
\end{equation}
where
 $B_{2k}(p)=\frac{(-1)^k}{(2k)!}\,(R_{p})^k,\
 B_{2k+1}(p)=\frac{(-1)^k}{(2k+1)!}\,(R_{p})^k A_p$
for $p\in M$.
\end{theorem}

The formulae (\ref{eq:main}) for few initial values of $k$,
$k=1,\ldots 3$, read as follows:
\begin{eqnarray}
\label{eq61}
 \int_M \sigma_1( A_p)\,{\rm d}V_F &=& 0,\\
\label{eq62}
 \int_M \big( \sigma_2(A_p) - \frac12\tr R_p\big)\,{\rm d}V_F &=& 0,\\
\label{eq63} \int_M
 \big(\sigma_3(A_p)-\frac{1}{2}\tr(A_p)\tr R_p +\frac{1}{3}\tr(R_p A_p)\big)\,{\rm d}V_F &=& 0.
\end{eqnarray}
The formulae (\ref{eq61}) and (\ref{eq62}) are well known for arbitrary foliated Riemannian manifolds, see the Introduction.
For $m=1$, (\ref{eq62}) reduces to the integral
of flag (Gauss) curvature, $\int_M K\,{\rm d}V_F=0$.

\begin{remark}\label{sec:nonsym}\rm
1. The compactness of $M$ in Theorem~\ref{thm:main} can be replaced by weaker conditions:
 $M$~is positively complete of finite $F$-volume, and has `bounded geometry' in the following sense:
\begin{equation}\label{E-bounded}
 \sup\nolimits_{\,p\in M} \|R_p\|_F<\infty,\quad \sup\nolimits_{\,p\in M} \|A_p\|_F<\infty.
\end{equation}
\ 2. Similar formulae exist for codimension-one foliations of on arbitrary
(non-locally symmetric with respect to $F$) Finsler manifolds.
They are more complicated since they contain terms which depend on covariant derivatives of $ R_p$.
More precisely, they contain just terms of the form
$ R_p^{(k)}$, where $ R_p^{(1)}=D_{\nu(p)} R_p$,
$ R_p^{(2)}=D_{\nu(p)}D_{\nu(p)}\, R_p$ and so on.
For the $F$-Jacobi tensor ${\bf Y}(t)$ we~get
\begin{equation*}
 {\bf Y}(t)=I_{m}+t A_p-\frac{t^2}{2!} R_p-\frac{t^3}{3!}( R_p A_p+ R_p^{(1)})
 +\frac{t^4}{4!}( R_p^2- R_p^{(2)}-2 R_p^{(1)} A_p)+\dots
\end{equation*}
The $t^3$ term of (\ref{eq5x}) becomes, compare (\ref{eq63}),
\begin{equation*}
 \int_M\big(\sigma_3( A_p)-\frac{1}{2}{\tr}( R_p)\tr( A_p)
 +\frac{1}{3}\tr( R_p\, A_p)-\frac{1}{6}\tr R_p^{(1)}\big)\,{\rm d}V_F=0.
\end{equation*}
In general, the $t^k$ term in (\ref{eq5x}) contains $ R_p^{(j)}$'s with $j\le k-2$.
\end{remark}

\begin{corollary}\label{cor-rank1}
Let $\calf$ be a codimension-one foliation on a $F$-locally symmetric
complete Finsler manifold $(M,F)$ of finite $F$-volume and bounded (in the sense of  \eqref{E-bounded}) geometry.
If ${\rm rank}(A_p)\le1$ for all $p\in M$
$($for example, $\calf$ is $F$-totally geodesic$)$ then the Riemannian curvature $ R_p$
vanishes identically provided that
$M$ has everywhere non-negative $($or, non-positive$)$ Ricci curvature $\Ric_p=\tr R_p$.
\end{corollary}

\proof Since in this case $\sigma_2(A_p)=0$, integral formula (\ref{eq62}) implies the claim.
\qed

\smallskip

Given a unit normal $\nu$ to $\calf$, denote by $Q_R$ the symmetric $(0,2)$-tensor in the rhs of \eqref{E-dt-Rm-Rm}.
Then,  see \eqref{E-dt-Ric-Ric},
\[
 \tr Q_R = I_{\nu}(\nabla^2_{\nu,\nu} \nu +C_\nu^\sharp(\nabla_{\nu}\,\nu))
 +2(\nabla_{\nu}I_{\nu})(\nabla_{\nu}\,\nu) -\tr\big(C^\sharp_{\,\nu}( C_\nu^\sharp + 2\,A^g)\big).
\]
 Define the 1-form $\theta_g$ by the equality
\[
 \theta_g(X)=g([X,\nu],\nu)\qquad (X\in TM).
\]
Note that $\nabla_\nu\,\nu=\theta_g^\sharp$ is the mean curvature of $\nu$-curves with respect to $g$.
 Comparing \eqref{eq:main} for $F$ and $g$, we obtain a series of integral formulas, the first two of which are given in the following.

\begin{theorem}
 Let $\tau(\nu)=\const$ on a codimension-one foliated Finsler space $(M,F)$. Then
\begin{eqnarray}\label{E-Q61}
 \int_M I_\nu(\nabla_\nu\,\nu)\,{\rm d}V_F \eq 0,\\
\label{E-Q62}
 \int_M \big( \sigma_2(C^\sharp_{\,\nu}) +(\tr A^g)(\tr C^\sharp_{\,\nu})-\tr(A^g C^\sharp_{\,\nu})
 - \frac12\tr Q_R\big)\,{\rm d}V_F \eq 0.
\end{eqnarray}
\end{theorem}

\proof
By \eqref{E-dt-A2}$_1$, $A = A^g+C^\sharp_{\,\nu}$, where $A=A_p$. Thus,
\eqref{E-Q61} follows from \eqref{eq61}, using \eqref{E-dt-A2}$_2$ and~Theorem~\ref{T-unit-n}.
 Note that by \eqref{eq-sigma-k1} with $k=1$ and \eqref{E-sigma22} (of Appendix),
 and by \eqref{E-dt-Ric-Ric}, we have
\begin{eqnarray*}
 \sigma_{2}(A_p) \eq \sigma_{2}(A^g) + \tr(A^g)\tr C^\sharp_{\,\nu}-\tr(A^g C^\sharp_{\,\nu}),\\
 \Ric_{\nu} \eq \tr R_p = \Ric^g_{\nu} +\tr Q_R.
\end{eqnarray*}
Thus, \eqref{E-Q62} follows from \eqref{eq62}, using \eqref{E-dt-A2}$_2$ and
\eqref{E-sigma22} with $k=2$ (of Appendix).
\qed

\subsection{Examples}
\label{sec:appl}

\textbf{Finsler manifolds of constant flag curvature}.
We provide examples, these of $(M,F)$ with constant flag curvature $K(\nu,P)$ on $M$,
i.e., such that  $R_p=K\,I_m$ for some $K\in\RR$.

a) For $(M,F)$ with zero flag curvature, $ R_p=0$, and we
obtain the Jacobi tensor of a simple form, linear in $t$:
${\bf Y}(t)=I_{m}+t A_p\ (t\ge0)$. Then (\ref{eq5x}) reduces to
 ${\rm Vol}_F(M)=\int_M\det(I_{m}+t A_p)\,{\rm d}V_F$.
 From this we obtain
 the Finsler generalization of the case $K=0$ of \cite[Theorem 1.1]{blr}, i.e.,
\begin{equation}\label{eq5e}
 \int_M \sigma_k( A_p)\,{\rm d}V_F=0,\quad k>0.
\end{equation}

b) Assume now that the flag curvature $K(\nu,P)$ of $(M, F)$ is constant and positive, say $K=1$.
Then $ R_p=I_m$ and one can rewrite the Taylor series for ${\bf Y}(t)\ (t\ge0)$ in the form
 ${\bf Y}(t)=\cos t\,\big(I_{m}+(\tan t) A_p\big)$.
Change of Variables Theorem for integration implies that the equality
\begin{equation*}
 {\rm Vol}_F(M) =(\cos t)^{m}\int_M\det\big(I_{m}+(\tan t) A_p\big)\,{\rm d}V_F
\end{equation*}
holds for arbitrary $t\ge0$ small enough. One can use the substitution $\tan t\to\tilde t$ and the
identity $\cos^2t=(1+\tilde t^2)^{-1}$ for further derivations.

c) The case of negative constant flag curvature $K(\nu,P)$ of $M$, say $K=-1$, is similar to the case~(b).
One can use the substitution $\tanh(t)\to\tilde t$ and the identity $\cosh^2t=(1-\tilde t^2)^{-1}$ for derivations.

The above yields the following extension of Theorem 1.1 in \cite{blr}.

\begin{corollary}\label{thm:003}
Let ${\cal F}$ be a transversally oriented codimension-one foliation on a Finsler mani\-fold $(M^{m+1},F)$
of finite $F$-volume and $\sup\nolimits_{\,p\in M} \|A_p\|_F<\infty$ (e.g. closed)
with a unit normal $\nu$ and condition $R_p =K I_m$. Then, for any $0\le k\le m$,
\begin{equation}\label{eq5f-b}
 \int_M \sigma_k( A_p)\,{\rm d}V_F = \left\{
 \begin{array}{cc}
 K^{k/2} \genfrac{(}{)}{0pt}{1}{\,m/2\,}{k/2}
 \,{\rm Vol}_F(M), & m,\,k \ {\rm even}, \\
 0, & m \ {\rm or}\ k \ {\rm odd}.
 \end{array}\right.
\end{equation}
\end{corollary}

\begin{remark}\label{Rem:nonsym}\rm
By Theorem~8.2.4 in \cite{mo}, if a Finsler manifold $M$ is closed
and has constant negative curvature then it is Randers.
\end{remark}

If $(M,F)$ is $F$-locally symmetric and the leaves of $\cal F$ are
$F$-\textit{totally geodesic} (i.e., $A_p = 0$)~then
\[
 {\bf Y}^{(2k+1)}(0)=0,\quad
 {\bf Y}^{(2k)}(0)=(- R_p)^k.
\]
Finally we get the $F$-Jacobi tensor
 ${\bf Y}(t)=I_{m}-\frac{t^2}{2!} R_p+\frac{t^4}{4!} R_p^2-\frac{t^6}{6!} R_p^3+\ldots$,
and (\ref{eq:main}) reduces to
\begin{equation*}
 \int_M \sum\nolimits_{\,\|{\lambda}\|=k}\sigma_{{\lambda}}\big({-}\frac{1}{2!}\,R_p,\ \frac1{4!}\,R_p^2\,,\ldots,\ \frac{(-1)^k}{(2k)!} R_p^k\big)\,{\rm d}V_F=0.
\end{equation*}
For codimension-one $F$-totally geodesic foliations on arbitrary  positively complete (or closed) Finsler manifolds of finite $F$-volume,
we get
\begin{eqnarray}\label{eq73}
 \nonumber
 \int_M \tr R_p\ {\rm d}V_F=0,\quad
 \int_M \tr R_p^{(1)}\,{\rm d}V_F=0,\\
 \int_M\big(\sigma_{2}( R_p)+\frac{1}{6}\,\tr R_p^2-\frac{1}{6}\,\tr R_p^{(2)}\big)\,{\rm d}V_F=0,
\end{eqnarray}
and so on. Equalities (\ref{eq73}) imply directly the following statement (see also Corollary~\ref{cor-rank1}).

\begin{corollary}\label{cor1}
Let $\calf$ be a codimension-one $F$-totally geodesic foliation on a $F$-locally symmetric po\-sitively complete Finsler manifold $(M,F)$ of finite $F$-volume and with condition \eqref{E-bounded}$_1$. Then $ R_p$ vanishes identically
provided that either $M$ has everywhere non-negative (or, non-positive) Ricci curvature $\Ric$, or $\sigma_2(R_p)$ is non-negative.
\end{corollary}

It has been observed in \cite{lw} that codimension-one foliations of compact negatively-Ricci
curved Riemannian spaces are far (in a sense) from being totally umbilical.
In~the case of an $F$-\textit{totally umbilical foliation}, $ A_p=H\,I_{m}$, therefore on a locally symmetric
Finsler space $(M,F)$ the following  can be derived from
(\ref{eq62})\,--\,(\ref{eq63})  etc. with the use of Lemma~\ref{lem1} of Appendix:
\begin{eqnarray} \label{eq75a}
 &&\int_M\big((m-1)(m-2) H^2 - \tr R_p\big)\,{\rm d}V_F = 0,\\
\label{eq75}
 &&\int_M  H\big(\frac{m(m-1)(m-2)}{3m-2} H^2 - \tr R_p\big)\,{\rm d}V_F = 0.
\end{eqnarray}
These integrals for $k$ even ((\ref{eq75a}), (\ref{eq75}), etc.)
contain polynomials depending on $ H^2$ only. If all the coefficients of such polynomials are positive, then
the polynomials are positive for all values of $ H$ and one may easily get
obstructions for existence of totally umbilical foliations on some Finsler manifolds.

\section{Codimension-one foliated Randers spaces}\label{sec:randers}

 Let $\calf$ be a transversally oriented codimension-one foliation of
 $M^{m+1}$ equipped with a Randers~metric
\[
 F(y)=\sqrt{a(y,y)} +\beta(y),\quad  \|\beta\,\|_\alpha<1,\quad \beta^\sharp\in\Gamma(T\calf).
\]
As before, let us write $a(\cdot, \cdot) = \langle\cdot, \cdot\rangle$.
Let $N$ be a unit $a$-normal vector field to $\calf$, i.e., $\<N,N\>=1$ and $\<N, v\>=0\ $ for $v\in T\calf$,
and $n$ an $F$-normal vector field to $\calf$ with the property $\<n,n\>=1$.
Denote by $\bar\nabla$ the Levi-Civita connection of the Riemannian metric $a$
and by $\nabla$ the Levi-Civita connection of the Riemannian metric $g=g_n$ on~$M$.
According to \cite[(1.15) and (1.19)]{cs} we have
\begin{eqnarray}\label{E-tauy}
 \tau(y) \eq (m+2)\log\sqrt{(1+\beta(y)/\alpha(y))\,c^{-2}}\,,\\
 \label{E-Iy}
 I_y(v) \eq \frac{m+2}{2 F(y)}\Big( \beta(v) -\frac{\<v,y\>\,\beta(y)}{\alpha^2(y)}\Big)\,.
\end{eqnarray}
In particular,
 $\tau(n)=0$ and $I_n(v) = \frac{m+2}{2\,c^4}\,\<\,\beta^\sharp -(c^2-1)\,n,\ v\>$.
Remark that for Randers spaces
\begin{equation*}
 C_n(u,v,w) = \frac1{m+2}\,\big( I_n(u)\,h_n(v,w)+I_n(v)\,h_n(u,w)+I_n(w)\,h_n(u,v)\big)\,,
\end{equation*}
where the angular form $h_n$ is given by
\begin{equation} \label{E-hy}
 h_n(u,v) = c^2\big( \<u,v\> -\<u,n\>\,\<v,n\>\big)\,,
\end{equation}
see \cite[(1.11) and (1.20)]{cs}.
Since $\sigma_F=c^{m+2}\sqrt{\det a_{ij}}$, see \cite[p.~6]{cs}, and
$\sqrt{\det g_{ij}(n)}=c^{m+2}\sqrt{\det a_{ij}}$, see \eqref{E-F001b},
the volume form of $F$ and canonical volume forms of Riemannian metrics $g$ and $a$ obey
\begin{equation}\label{E-F001vol}
 {\rm d}V_F=c^{m+2}{\rm d}V_a,\qquad
 {\rm d} V_g = c^{m+2}{\rm d} V_a ,\qquad
 {\rm d}V_F= {\rm d}V_g.
\end{equation}
Let $Z=\nabla_\nu\,\nu$ (which is dual of $\theta_g$ in Sect.~\ref{sec:IF}) and $\bar Z=\bar\nabla_N\,N$ be the curvature vectors of $\nu$-curves and $N$-curves for Riemannian metrics $g$ and $a$, respectively.

\subsection{The shape operators of $g$ and $a$}

The shape operators of $\calf$ with respect to the metrics $a$ and $g$ are defined as follows:
\[
 \bar A(u)=-\bar\nabla_u\,N,\quad A^g(u)=-\nabla_u\,\nu\,,
\]
where $u\in T\calf$ and $\nu =c^{-2}n = c^{-1}(N - c^{-1}\beta^\sharp)$ with $c=\sqrt{1-\|\beta\,\|^2_\alpha}>0$\,.

The derivative $\bar\nabla u:TM\to TM$ is defined by $(\bar\nabla u)\,(v)=\bar\nabla_v\,u = \bar\nabla_v\,u$,
where $v\in TM$. The~conjugate derivative $(\bar\nabla u)^t: TM\to TM$ is defined by
$\<(\bar\nabla\,u)^t(v),w\>=\<v,(\bar\nabla\,u)(w)\>$ for all $v,w\in TM$.
 The \textit{deformation tensor} $\overline{\rm Def}$,
\[
 2\,\overline{\rm Def}_u =\bar\nabla u+(\bar\nabla u)^t,
\]
measures the degree to which the flow of a vector field $u\in\Gamma(T M)$ distorts the metric $a$.
The~same notation $\overline{\rm Def}_u$ will be used for its dual (with respect to $a$) $(1,1)$-tensor.
Set $\overline{\rm Def}^\top_u(v)=(\,\overline{\rm Def}_u(v))^\top$.
For $\beta\ne0$, let
\[
 \bar A(\beta^\sharp)^{\bot\beta}=\bar A(\beta^\sharp)-\<\bar A(\beta^\sharp), \beta^\sharp\>\beta^\sharp\cdot\|\beta^\sharp\|^{-2}_\alpha
\]
be the projection of $\bar A(\beta^\sharp)$ on $(\beta^\sharp)^\bot$.
Note that $\lim_{\,\beta\to0}\bar A(\beta^\sharp)^{\bot\beta}=0$.

\begin{proposition}\label{L-Dx}
Let $\,\beta(N)=0$ on $M$. Then on $T\calf$ we have
\begin{equation}\label{E-A-bar-A}
 c\,A^g = \bar A - c^{-2}(c\,N-\beta^\sharp)(c)I_m +\,c^{-1}
 (\overline{\rm Def}_{\beta^\sharp})^\top_{|T\calf}
 + U_1^\flat\otimes\beta^\sharp +U_2\otimes\beta\,,
\end{equation}
where
\begin{eqnarray}\label{E-A-bar-AU}
\nonumber
 U_1\eq-\frac12\,c^{-2}\big( (c\,N-\beta^\sharp)(c)\,\beta^\sharp
 -2\,c^{-1}(\,\overline{\rm Def}_{\beta^\sharp}\,\beta^\sharp)^\top
 -\bar\nabla^\top_{N-c^{-1}\beta^\sharp}\beta^\sharp  \\
\nonumber
 && +\,c\,\bar Z +c\,\beta(\bar Z)\,\beta^\sharp -\bar A(\beta^\sharp)^{\bot\beta}\big), \\
 U_2\eq\frac12\,\big(\bar\nabla_{N-c^{-1}\beta^\sharp}^\top\,\beta^\sharp
 -c\,\bar Z -\bar A(\beta^\sharp)^{\bot\beta}\,\big)\,.
\end{eqnarray}
\end{proposition}

\proof
By the well-known formula for Levi-Civita connection of $g$,
using equalities $g(u,n)=0=g(v,n)$ and $g([u,v],n)=0$, we have
\begin{equation}\label{E-LC-g}
 2\,g(\nabla_u\,n, v) = n(g(u,v)) +g([u,n],v) +g([v,n],u)\quad
 (u,v\in T\calf).
\end{equation}
One may assume $\bar\nabla_X^\top\,u=\bar\nabla_X^\top\,v=0$ for all $X\in T_pM$ at a given point $p\in M$.
Using \eqref{E-gnv} with $u=[u,n]$ and $v=v$, we obtain
\begin{eqnarray*}
 n(g(u,v))\eq n(c^2(\<u,v\> -\beta(u)\,\beta(v)))\\
 \eq n(c^2)(\<u,v\>-\beta(u)\beta(v)) -c^2\beta(u)(\bar\nabla_n\,\beta)(v) -c^2(\bar\nabla_n\,\beta)(u)\beta(v),\\
 g([u,n],v) \eq c^2\big(\<[u,n],v\> +\beta(v)\<[u,n]),\,n\>\big) \\
 \eq -c^2\<c\,\bar A(u) +\bar\nabla_u\,\beta^\sharp, v\> +c^3\<\bar A(\beta^\sharp)+c\bar Z,\,u\>\,\beta(v),\\
 g([v,n],u) \eq c^2\big(\<[v,n],u\> +\beta(u)\<[v,n]),\,n\>\big)\\
 \eq -c^2\<c\,\bar A(v) +\bar\nabla_v\,\beta^\sharp, u\> +c^3\<\bar A(\beta^\sharp)+c\bar Z,\,v\>\,\beta(u).
\end{eqnarray*}
Substituting the above into \eqref{E-LC-g}, we find
\begin{eqnarray}\label{E-DDD}
\nonumber
 &&\hskip-10mm 2\,g(\nabla_u\,n, v) = n(c^2)\big(\<u,v\> -\beta(u)\beta(v)\big)
 -2\,c^3\<\bar A(u),v\> -2\,c^2\<\overline{\rm Def}_{\beta^\sharp}(u),v\> \\
\nonumber
 &&\hskip-10mm -c^2(\bar\nabla_n\,\beta)(u)\,\beta(v) -c^2\beta(u)(\bar\nabla_n\,\beta)(v)
 +c^3\<\bar A(\beta^\sharp)+c\bar Z,\,u\>\,\beta(v) \\
 && +\,c^3\beta(u)\<\bar A(\beta^\sharp)+c\bar Z,\,v\>.
\end{eqnarray}
 From \eqref{E-DDD}, assuming $g(\nabla_u\,n, v)=\<\mathfrak{D}(u),\,v\>$ and using Lemma~\ref{L-zZ}, we get
\begin{eqnarray}\label{E-Dx}
 -2\,c^4 A^g(u)\eq 2\,\mathfrak{D}(u)+c^{-2}\<2\,\mathfrak{D}(u),\,\beta^\sharp\>\,\beta^\sharp\,,
\end{eqnarray}
where $\mathfrak{D}: T\calf\to T\calf$ is a linear operator, and
\begin{eqnarray}\label{E-Dx2}
\nonumber
 2\,\mathfrak{D}(u) \eq n(c^2)\,(u - \beta(u)\beta^\sharp) -2\,c^3\bar A(u) -2\,c^2\,(\,\overline{\rm Def}_{\beta^\sharp}(u))^\top \\
\nonumber
 \minus c^2(\bar\nabla^\top_n\,\beta)(u)\,\beta^\sharp -c^2\beta(u)\bar\nabla_n^\top\,\beta^\sharp
 +c^3\<\bar A(\beta^\sharp)+c\bar Z,\,u\>\,\beta^\sharp \\
 \plus c^3\beta(u)(\bar A(\beta^\sharp) +c\bar Z).
\end{eqnarray}
From \eqref{E-Dx2} we get
\begin{eqnarray}\label{E-Dx3}
\nonumber
 2\,\<\mathfrak{D}(u),\beta^\sharp\> \eq n(c^2)\,c^2\beta(u)
 -2\,c^3\<\bar A(\beta^\sharp),u\> -2\,c^2\,\<\overline{\rm Def}_{\beta^\sharp}(\beta^\sharp),u\> \\
\nonumber
 \minus c^2(1-c^2)\,(\bar\nabla^\top_n\,\beta)(u) +c^3 n(c)\beta(u) +c^3(1-c^2)\,\<\bar A(\beta^\sharp)+c\bar Z,\,u\> \\
 \plus c^3\<\bar A(\beta^\sharp)+c\bar Z,\beta^\sharp\>\,\beta(u).
\end{eqnarray}
 From \eqref{E-Dx}\,--\,\eqref{E-Dx3} we obtain
\begin{eqnarray*}
 && c\,A^g = \bar A -c^{-1}\,(N-c^{-1}\,\beta^\sharp)(c)\,I_m
 c^{-1}\,(\overline{\rm Def}_{\beta^\sharp})^\top_{|T\calf}
 \\
 && -\,\frac12\,c^{-2}\,\big( (c\,N-\beta^\sharp)(c)\,\beta^\sharp -2\,c^{-1}(\,\overline{\rm Def}_{\beta^\sharp}\,\beta^\sharp)^\top
 \!-\bar\nabla^\top_{N-c^{-1}\beta^\sharp}\,\beta^\sharp
 +c\,\bar Z +c\,\<\bar Z,\,\beta^\sharp\>\,\beta^\sharp \\
 && -\,\bar A(\beta^\sharp) +\<\bar A(\beta^\sharp),\beta^\sharp\>\,\beta^\sharp
 \big)^\flat\otimes\beta^\sharp
 +\frac12\,\big(\bar\nabla_{N-c^{-1}\beta^\sharp}^\top\,\beta^\sharp -c\bar Z -\bar A(\beta^\sharp)\big)\otimes\beta .
\end{eqnarray*}
 From the above the expected \eqref{E-A-bar-A}\,--\,\eqref{E-A-bar-AU} follow.
\qed

\begin{corollary}\label{C-shape}
 Let $\,\beta(N)=0$. If $\|\beta\|_\alpha=\const$ then on $T\calf$ we have
\begin{eqnarray}\label{E-A-bar-A-c}
\nonumber
 && c\,A^g = \bar A +c^{-1}\,(\overline{\rm Def}_{\beta^\sharp})^\top_{|T\calf}
 +\frac12\,\big(\bar\nabla_{N-c^{-1}\beta^\sharp}^\top\,\beta^\sharp -c\bar Z
 -\bar A(\beta^\sharp)^{\bot\beta}\big)\otimes\beta \\
\nonumber
 && +\,\frac12\,c^{-2}\,\big( 2\,c^{-1}\overline{\rm Def}\,^\top_{\beta^\sharp}(\beta^\sharp)
 +\bar\nabla^\top_{N-c^{-1}\beta^\sharp}\,\beta^\sharp +\bar A(\beta^\sharp)^{\bot\beta} \\
 && -\,c\,\bar Z -c\,\<\bar Z,\,\beta^\sharp\>\,\beta^\sharp \big)^\flat\otimes\beta^\sharp.
\end{eqnarray}
 If, in particular, $\bar\nabla\beta=0$ $($i.e., $F$ is a Berwald structure$)$ then
\begin{equation}\label{E-A-bar-A-fin}
 c\,A^g = \bar A -\frac12\,\big(\bar A(\beta^\sharp)^{\bot\beta}+c\bar Z\big)\otimes\beta
 +\frac12\,c^{-2}\,\big( \bar A(\beta^\sharp)^{\bot\beta}
 -c\,\bar Z -c\,\<\bar Z,\,\beta^\sharp\>\,\beta^\sharp \big)^\flat\otimes\beta^\sharp.
\end{equation}
\end{corollary}

\subsection{The Riemann curvature of $g$ and $a$}

In this section we study relationship between Riemann curvature of two metrics, $g$ and $a$, on a~Randers space.

\begin{proposition}\label{P-ZZ} For a codimension-one foliation of $M$ with Riemannian metrics $g$ and $a$ we have
\begin{eqnarray}\label{E-ZZ}
 Z \eq c^{-2}\bar Z -c^{-3}\,\bar\nabla^\top c +c^{-4}\,\beta(\bar Z-c^{-1}\,\bar\nabla^\top c) \,\beta^\sharp,\\
\label{E-Csharp}
 C^\sharp_n \eq c^{-2}\bar C +c^{-4}(\beta\circ\bar C)\otimes\beta^\sharp,
\end{eqnarray}
where
\begin{eqnarray}
\nonumber
 2\,\bar C \eq {\rm Sym}(\beta\otimes\bar Z)
 +c^{-3}\big(c\,\beta(\bar Z) -2\,\beta^\sharp(c) -n(c)\big)\big(I_m -\beta\otimes\beta^\sharp\big)\\
 \nonumber
 \minus c^{-1}{\rm Sym}(\beta\otimes\bar\nabla^\top c) +c^{-1}(\beta^\sharp(c) +n(c))\big(I_m -3\,\beta\otimes\beta^\sharp\big).
\end{eqnarray}
We also have
\begin{eqnarray}\label{E-nablaZZ}
 &&\quad \<\bar\nabla_u \bar Z,v\> = \<\bar\nabla_v\bar Z,u\>,\quad
 g(\nabla_u Z,v) = g(\nabla_v Z,u)\quad (u,v\in T\calf),\\
\label{E-RnRnu-R}
 &&\ \ \bar R_N = (\overline{\rm Def}_{\bar Z})^\top_{|T\calf} {+}\bar\nabla_N \bar A {-}\bar A^2
 {-}\bar Z^\flat\otimes\bar Z,\quad
 R^g_\nu = ({\rm Def}_Z)^\top_{|T\calf}
 {+}\nabla_\nu A {-}A^2 {-}Z^\flat\otimes Z.
\end{eqnarray}
\end{proposition}

\proof
Extend $X\in T_p\calf$ at a point $p\in M$ onto a neighborhood of $p$ with the property $(\bar\nabla_Y\,X)^\top=0$
for any $Y\in T_pM$. By the well known formula for the Levi-Civita connection, we obtain at $p$:
\[
 g(Z,X) = g([X,\nu],\nu).
\]
Then, using the equalities $\nu=c^{-1}N-c^{-2}\beta^\sharp$ and $[X, fY]=X(f)Y+f[X, Y]$, we calculate
\begin{eqnarray*}
 g([X,\nu],\nu) \eq c^{-4}X(c)\,g(N,\beta^\sharp) -c^{-3}X(c)\,g(N,N) \\
 \plus c^{-2} g([X,N],\,N) - c^{-3} g([X,N], \beta^\sharp).
\end{eqnarray*}
Note that
\[
 [X,N] = \bar\nabla_X N -\bar\nabla_N X =-\bar A(X) -\<\bar\nabla_N X,\,N\>\,N
 = -\bar A(X) +\<\bar Z,\,X\>\,N
\]
and $N=c\nu +c^{-1}\beta^\sharp$. Then, by Lemma~\ref{L-c-value} and the equalities
\begin{eqnarray*}
 g(\beta^\sharp,\beta^\sharp) \eq c^2(\<\beta^\sharp,\,\beta^\sharp\> -\beta(\beta^\sharp)^2) = c^4(1-c^2),\\
 g(N,\beta^\sharp) \eq g(c\nu+c^{-1}\beta^\sharp,\,\beta^\sharp) =c^{-1}g(\beta^\sharp,\beta^\sharp)= c^3(1-c^2),\\
 g(N,N) \eq g(c\nu+c^{-1}\beta^\sharp,\,c\nu+c^{-1}\beta^\sharp) = c^2 +c^{-2}g(\beta^\sharp,\beta^\sharp) = c^2(2-c^2),
\end{eqnarray*}
we obtain
\begin{eqnarray*}
 g([X,N],\,N) \eq -\<\bar A(\beta^\sharp),\,X\> +\<\bar Z, X\>\,g(N,N)
 = c^2\<(2-c^2)\bar Z -c\bar A(\beta^\sharp),\ X\>,\\
 g([X,N], \beta^\sharp) \eq -\<\bar A(\beta^\sharp),\,X\> +\<\bar Z, X\>\,g(N,\beta^\sharp)
 = c^3\<(1-c^2)\bar Z -c\bar A(\beta^\sharp),\ X\>.
\end{eqnarray*}
Hence,
\begin{eqnarray*}
 g(Z,X) \eq  -c^{-1}X(c) +\<\bar Z, X\> = \<\bar Z -c^{-1}\bar\nabla c,\ X\>.
\end{eqnarray*}
By Lemma~\ref{L-zZ}, we get \eqref{E-ZZ}.
From \eqref{E-Iy}\,--\,\eqref{E-hy}, \eqref{E-ZZ} and a bit of help from Maple program we find
\begin{eqnarray*}
 2\,C_n(u,v, Z) \eq \<\bar Z,u\>\,\beta(v)+\<\bar Z,v\>\,\beta(u) \\
 \plus c^{-3}(c\,\beta(\bar Z) -2\,\beta^\sharp(c) -n(c))\big(\<u,v\> -\beta(u)\beta(v)\big)\\
 \minus c^{-1}(u(c)\,\beta(v) +v(c)\,\beta(u))
 +c^{-1}(\beta^\sharp(c) +n(c))\big(\<u,v\> - 3\,\beta(u)\,\beta(v)\big).
\end{eqnarray*}
Using $g(C^\sharp_n(u),v)=\<\bar C(u),v\>$, where $C^\sharp_n$ is $g$-dual to $C_n(\cdot,\cdot, \nabla_n\,n)$, and
\begin{eqnarray*}
 2\,\bar C(u) \eq \<\bar Z,u\>\,\beta^\sharp+\beta(u)\bar Z
 +c^{-3}(c\,\beta(\bar Z) -2\,\beta^\sharp(c) -n(c))\big(u -\beta(u)\beta^\sharp\big)\\
 \minus c^{-1}(u(c)\,\beta^\sharp +\beta(u)\,\bar\nabla^\top c)
 +c^{-1}(\beta^\sharp(c) +n(c))\big(u -3\,\beta(u)\,\beta^\sharp\big),
\end{eqnarray*}
we apply Lemma~\ref{L-zZ} to get \eqref{E-Csharp}.

We shall prove \eqref{E-nablaZZ} and \eqref{E-RnRnu-R} for $a$.
It is sufficient to show that
\begin{equation}\label{E-RuNvN-2}
 \<\bar R(u,N)N,v\> = \<(\bar\nabla_N\,\bar A - \bar A^2)(u),v\> - \<\bar Z, u\>\<\bar Z, v\> +\<\bar\nabla_u\,\bar Z, v\>,\quad
 u,v\in T\calf.
\end{equation}
Since the left hand side of \eqref{E-RuNvN-2} is symmetric, we obtain
$\<\bar\nabla_u\bar Z,v\>=\<\bar\nabla_v \bar Z,u\>$, see \eqref{E-RnRnu-R}$_1$ and \eqref{E-nablaZZ}$_1$.
Indeed,
\begin{eqnarray*}
 &&\<\bar R(u,N)N,v\> = \<\bar\nabla_u\bar\nabla_N N,v\> -\<\bar\nabla_N\bar\nabla_u N,v\>
 -\<\bar\nabla_{\bar\nabla_u N -\bar\nabla_N u} N,v\>\\
 && = \<\bar\nabla_u \bar Z,v\> +\<\bar\nabla_N(\bar A(u)),v\> -\<\bar A^{\,2}(u),v\>
 +\<\bar\nabla_{\<\bar\nabla_N\,u,N\>N}\,N,v\> -\<\bar A(\bar\nabla^\top_N\,u),v\>\\
 && = \<(\bar\nabla_N\,\bar A - \bar A^{\,2})(u),v\> - \<\bar Z, u\>\<\bar Z, v\> +\<\bar\nabla_u\,\bar Z, v\>,
\end{eqnarray*}
that completes the proof of \eqref{E-RuNvN-2}.
The proof of \eqref{E-nablaZZ}$_2$  and \eqref{E-RnRnu-R}$_2$ (for the metric $g$)  is similar.
\qed

\smallskip
By \eqref{E-Csharp}, the~equality $C^\sharp_n=0$ is independent of the condition $\bar\nabla\beta=0$.
Moreover, we have the following.

\begin{corollary}\label{C-sharp-n}\rm
Let $m>3$ and $c=\const$. Then $C^\sharp_n=0$ if and only if $\bar Z=0$.
\end{corollary}

\proof By our assumptions,
$\bar C = \frac12\,{\rm Sym}(\beta\otimes\bar Z)+\frac12\,c^{-2}\,\beta(\bar Z) \big(I_m -\beta\otimes\beta^\sharp\big)$.
Hence, $C^\sharp_n=0$ reads
\[
 \beta(\bar Z)I_m =\beta(\bar Z)\,\beta\otimes\beta^\sharp
 -c^2\,{\rm Sym}(\beta\otimes\bar Z) -2\,(\beta\circ\bar C) \otimes\beta^\sharp.
\]
Since the matrix $\beta(\bar Z)I_m$ is conformal,
while the matrix in the right hand side of above equality has the form $\omega\otimes\beta^{\sharp} -c^2\bar Z^{\,\bot\beta}\otimes\beta$
and rank $\le 3$, for $m>3$ we obtain
\[
 \beta(\bar Z)=0,\quad
 {\rm Sym}(\beta\otimes\bar Z) +2\,c^{-2}(\beta\circ\bar C)\otimes\beta^\sharp =0.
\]
By the first condition, $\bar Z\perp\beta^\sharp$; thus, the second condition yields
$\bar Z=0$  (that is, $\calf$ is a Riemannian foliation for the metric $a$) and $\bar C=0$.
The converse claim follows directly from \eqref{E-Csharp} and the definition of $\bar C$.
\qed

\begin{remark}\label{C-3beta}\rm
In \cite{cs2,sh2} one may find coordinate presentations of
$R_y$ through $\bar R_y$ for all $y\in TM$.
For example, if $\bar\nabla\beta=0$ $($i.e., $F$ is a Berwald structure$)$ then $R_y(u) = \bar R_y(u)$ for all $u$.
Alternative formulas with relationship between $R_\nu$ and $\bar R_\nu$ follow from \eqref{E-RnRnu-R}, where
$A^g$ and $Z$ are expressed using $\bar A$ and $\bar Z$ given in Propositions~\ref{L-Dx} and \ref{P-ZZ}.
\end{remark}

\subsection{Around the Reeb and Brito-Langevin-Rosenberg formula}

In results  of this section, a closed manifold can be replaced by a complete manifold of finite volume
with bounded geometry, see conditions \eqref{E-bounded}.
Based on \eqref{eq:main} and \eqref{eq5f-b}, one may produce a sequence of similar formulae for Randers spaces.
We~will discuss first two of them (i.e., $k=1,2$).

\begin{remark}\rm
In \cite{re}, G.\,Reeb proved that the total mean curvature of the leaves of a~codimen\-sion-one foliation on a closed Riemannian manifold equals zero.
Note that $\tr\overline{\rm Def}\,^\top_{\beta^\sharp}=\overline{\Div}\,\beta^\sharp +\beta(\bar Z)$,
where $\bar Z=\bar\nabla_N\,N$  is the curvature vector of $N$-curves for the metric $a$.
Using notations of Appendix, we find from \eqref{E-A-bar-AU},
\[
 \beta(U_1)=-\frac{2-c^2}{2\,c}\,N(c) -\frac{1}{2}\,\beta^\sharp(c) -\frac{2-c^2}{2\,c}\,\beta(\bar Z),\quad
 \beta(U_2)=-\frac{1}{2}\,(c\,N -\beta^\sharp)(c)-\frac{1}{2}\,c\,\beta(\bar Z).
\]
Hence,
\[
 \beta(U_1)+\beta(U_2)=-c^{-1}(N(c)+\beta(\bar Z)).
\]
Tracing \eqref{E-A-bar-A}, we get
\begin{eqnarray*}
  c\,\sigma_1(A^g) = \sigma_1(\bar A) -(m+1)\,c^{-1}N(c) +m\,c^{-2}\beta^\sharp(c) +c^{-1}\,\overline{\Div}\,\beta^\sharp.
\end{eqnarray*}
The volume forms of $g$ and $a$ obey ${\rm d}V_g=c^{m+2}\,{\rm d}V_a$, see \eqref{E-F001vol}.
Using the Reeb formula for metric $g$,
\[
\int_M \sigma_1(A^g)\,{\rm d}V_g=0,
\]
the equality $\overline\Div(c^{m}\beta^\sharp)=c^{m}\,\overline\Div\,\beta^\sharp +\beta^\sharp(c^m)$
and the Divergence Theorem, we~get
\begin{equation}\label{E-IF1-Randers0}
 \int_M \big(c^{m+1}\sigma_1(\bar A) -N(c^{m+1})\big)\,{\rm d}V_a =0\,,
\end{equation}
which for $\beta=0$ is the Reeb formula for metric $a$.
Remark that \eqref{E-IF1-Randers0} is a particular case of a general formula
for any $f\in C^2(M)$, see \cite[Lemma~2.5]{rw1}:
\[
 \int_M(f\,\sigma_1(\bar A)-N(f))\,{\rm d}V_a = 0.
\]
\end{remark}

The next results concern Brito-Langevin-Rosenberg type formulas for foliated Randers spaces.

The \textit{Newton transformations} $T_k(A)\ (0\le k\le m)$ of an $m\times m$ matrix $A$ (see \cite{rw1}) are defined either inductively by
$T_0(A)=I_m,\ T_k(A)=\sigma_k(A) I_m-A\,T_{k-1}(A)\ (k\ge1)$ or explicitly~as
\begin{equation*}
 T_k(A) = \sigma_k(A) I_m-\sigma_{k-1}(A)\,A+\ldots+(-1)^k A^k,\quad 0\le k\le m,
\end{equation*}
and we have $T_k(\lambda\,A)=\lambda^k\,T_k(A)$ for $\lambda\ne0$.
Observe that if~a rank-one matrix $A:=U\otimes\beta$ (and similarly for $A:=\omega\otimes\beta^\sharp$) has zero trace,
i.e., $\beta(U)=0$, then
\[
 A^2=U(\beta^\sharp)^t \cdot U(\beta^\sharp)^t = U\beta(U)\,(\beta^\sharp)^t=\beta(U)\,A=0.
\]
Note that for $c=\const$ we have, see \eqref{E-Csharp},
$C^\sharp_n = c^{-2}\bar C +c^{-4}(\beta\circ\bar C)\otimes\beta^\sharp$, where
$C^\sharp_n =c^2 C^\sharp_\nu$ and
\begin{equation*}
 2\,\bar C = {\rm Sym}(\beta\otimes\bar Z) +c^{-2}\beta(\bar Z)(I_m -\beta\otimes\beta^\sharp).
\end{equation*}

\begin{theorem}\label{T-Randers-Kconst}
Let $(M^{m+1},\alpha+\beta)$ be a codimension-one foliated closed Randers space
with constant sectional curvature $\bar K$ of $a$.
If a nonzero vector field $\beta^\sharp\in\Gamma(T\calf)$ obeys $\,\bar\nabla\beta=0$, then $\bar K=0$ and
for $1\le k\le m$ we have
\begin{eqnarray}\label{E-IF-Randers-k}
\nonumber
 &&\hskip-12mm\int_M \Big(\sum\nolimits_{j>0}\sigma_{k-j,j}(\bar A, c\,C^\sharp_\nu)
 +\<T_{k-1}(\bar A+c\,C^\sharp_\nu)(\beta^\sharp),\, U_1\> \\
 && +\,\big\<T_{k-1}(\bar A+c\,C^\sharp_\nu+U_1^\flat\otimes\beta^\sharp)(U_2),\,\beta^\sharp\big\>
 \Big)\,{\rm d}V_a = 0,
\end{eqnarray}
where
 $U_1=\frac12\,c^{-2}(\bar A(\beta^\sharp) -c\bar Z)$,
 $U_2=-\frac12(\bar A(\beta^\sharp) +c\bar Z)$.
Moreover, if $m>3$ and $\bar Z=0$ then
\begin{equation}\label{eq5f-b-R2}
 \int_M \big\< \big(c^{-2}T_{k-1}(\bar A) -T_{k-1}(\bar A +\frac12\,c^{-2}\bar A(\beta^\sharp)^\flat
 \otimes\beta^\sharp)\big)(\bar A(\beta^\sharp)),\ \beta^\sharp\big>\,{\rm d}V_a =0.
\end{equation}
\end{theorem}

\proof
By our assumptions, $c={\rm const}$ and $\bar R(x,y)z=\bar K(\,\<y,z\>\,x-\<x,z\>\,y\,)$.
Hence, on $T\calf$
\[
 \bar R_N=\bar K I_m,\quad
 \bar R_{\beta^\sharp}=(1-c^2)\bar K I_m,\quad
 \bar R(\cdot,N)\beta^\sharp=0.
\]
If $\,\bar\nabla\beta=0$ then $\bar R(U,\beta^\sharp,\beta^\sharp,U) = 0$ and $\bar K(U\wedge\beta^\sharp) = 0$
for all $U\perp\beta^\sharp$; hence, in our case, $\bar K=0$.
 By Remark~\ref{C-3beta}, $R_y = \bar R_y$ for all $y\in TM_0$; hence, $R_y=0$.
Since $\bar\nabla\beta^\sharp=0$, we obtain $\beta(\bar Z)=0$ and $\<\bar A(\beta^\sharp),\beta^\sharp\>=0$:
\begin{eqnarray*}
 && \<\beta^\sharp,\,\bar Z\> = \<\beta^\sharp,\,\bar\nabla_N N\> = -\<\bar\nabla_N\,\beta^\sharp,\,N\>=0,\\
 && \<\bar A(\beta^\sharp),\beta^\sharp\>
 =-\<\beta^\sharp,\,\bar\nabla_{\beta^\sharp}\,N\>
 =\<\bar\nabla_{\beta^\sharp}\,\beta^\sharp,\,N\>=0.
\end{eqnarray*}
By \eqref{E-dt-A2} and Corollary~\ref{C-shape},
\[
 c\,A = c\,A^g+c\,C^\sharp_\nu = \bar A+c\,C^\sharp_\nu+A_1+A_2,
\]
where $A_1 = U_1^\flat\otimes\beta^\sharp$ and $A_2 = U_2\otimes\beta$ are rank $\le 1$ matrices
(since $\<U_i,\beta^\sharp\>=0$).
By~Corollary~\ref{C-AplusB} of Appendix, we have
\begin{eqnarray}\label{E-2sigma-k}
\nonumber
 c^k\sigma_k(A) \eq \sigma_k(\bar A) +\sum\nolimits_{j>0}\sigma_{k-j,j}(\bar A, c\,C^\sharp_\nu)
 +U_1\big(T_{k-1}(\bar A+c\,C^\sharp_\nu)(\beta^\sharp)\big) \\
 \plus \beta\big(T_{k-1}(\bar A+c\,C^\sharp_\nu+A_1)(U_2)\big).
\end{eqnarray}
Recall that  ${\rm d}V_F=c^{m+2}\,{\rm d}V_a$, see \eqref{E-F001vol}.
Comparing \eqref{eq5f-b} (when $K=0$)~with
\begin{equation*}
 \int_M \sigma_k(\bar A_p)\,{\rm d}V_a = 0,
\end{equation*}
we find \eqref{E-IF-Randers-k}.
By Corollary~\ref{C-sharp-n}, if $m>3$, $\bar Z=0$ then $C^\sharp_\nu=0$;
hence, \eqref{E-IF-Randers-k} yields~\eqref{eq5f-b-R2}.
\qed

\begin{remark}\label{Ex-k1}\rm
For $k=1$, \eqref{E-IF-Randers-k} yields the Reeb type formula
\begin{equation*}
 \int_M
 \sigma_1(C^\sharp_\nu)
 \,{\rm d}V_a = 0.
\end{equation*}
\end{remark}

\begin{corollary}
Let $(M^{m+1},\alpha+\beta),\ m>3$, be a codimension-one foliated closed Randers space
with constant sectional curvature~$\bar K$ of $a$. If $\bar Z=0$ and a
nonzero vector field $\beta^\sharp\in\Gamma(T\calf)$ obeys $\bar\nabla\beta=0$
then $\bar K=0$ and
$\bar A(\beta^\sharp)=0$
at any point of $M$.
If, in addition, $\calf$ is totally umbilical $(\bar A = \bar H\cdot I_m)$ then $\calf$ is totally geodesic.
\end{corollary}

\proof For $k=2$, the integrand in \eqref{eq5f-b-R2} reduces to
 $\frac{c^2-1}{4\,c^{2}}\,\|\bar A(\beta^\sharp)\|^2$. Thus, when $c\ne1$, the claim follows.

Nevertheless, we will give alternative proof with use of integral formula \eqref{eq62}.
Our Randers space $(M,\alpha+\beta)$ is now Berwald.
For the rank 1 matrices $A_1=U_1^\flat\otimes\beta^\sharp$ and $A_2=U_2\otimes\beta$,
where $U_1=\frac12\,c^{-2}\bar A(\beta^\sharp)$ and $U_2=-\frac12\bar A(\beta^\sharp)$
and $\<\bar A(\beta^\sharp),\beta^\sharp\>=0$,
see \eqref{E-A-bar-A-fin} with $\bar Z=0$, we~have
\begin{eqnarray*}
 &&\tr(A_1 A_2)=\<U_1,U_2\>\,\beta(\beta^\sharp)=\frac{c^2-1}{4\,c^{2}}\,\|\bar A(\beta^\sharp)\|^2_\alpha,\\
 &&\tr(\bar A A_1)=\<U_1,\bar A(\beta^\sharp)\>=\frac{1}{2\,c^{2}}\,\|\bar A(\beta^\sharp)\|^2_\alpha,\\
 &&\tr(\bar A A_2) =\<U_2,\bar A(\beta^\sharp)\>=-\frac{1}{2}\,\|\bar A(\beta^\sharp)\|^2_\alpha.
\end{eqnarray*}
Thus, $\tr(A_1 A_2+\bar A A_1+\bar A A_2)=\frac{1-c^2}{4\,c^{2}}\|\bar A(\beta^\sharp)\|^2$.
By the identity for square matrices
\begin{eqnarray*}
 \sigma_2(\sum\nolimits_{\,i} A_i) \eq \frac12\,\tr^2(\sum\nolimits_{\,i} A_i)-\frac12\,\tr((\sum\nolimits_{\,i} A_i)^2)\\
 \eq\sum\nolimits_{\,i}\sigma_2(A_i)+\sum\nolimits_{\,i<j} \big((\tr A_i)(\tr A_j) -\tr(A_i A_j)\big)\,,
\end{eqnarray*}
and $\sigma_2(A_1)=\sigma_2(A_2)=0$, by~the above and since
$c\,A=c\,A^g=\bar A+A_1+A_2$, we get
\begin{eqnarray*}
 c^2\sigma_2(A)\eq c^2\sigma_2(A^g) =\sigma_2(\bar A) +\frac14\,(c^{-2}-1)\,\|\bar A(\beta^\sharp)\|^2_\alpha.
\end{eqnarray*}
 From the integral formulae, \eqref{eq5e}, for $F$ and for Riemannian metric $a$,
\[
 \int_M\sigma_2(\bar A)\,{\rm d}V_a=0,\quad
 \int_M\sigma_2(A)\,{\rm d}V_F=0,
\]
where the volume forms are related by ${\rm d}V_F=c^{m+2}{\rm d}V_a$, see \eqref{E-F001b}, we find that
\[
 (c^{-2}-1)\int_M \|\bar A(\beta^\sharp)\|^2_\alpha\,{\rm d}V_a = 0.
\]
Since $c\ne1$ (for $\beta\ne0$), we obtain $\bar A(\beta^\sharp) =0$.
\qed

\smallskip

Similar integral formulae exist for codimension one totally umbilical
(i.e., $\bar A = \bar H I_m$, where $\bar H=\frac1m\tr\bar A$) and totally geodesic foliations.
Notice that non-flat closed Riemannian manifolds of constant curvature do not admit such foliations.

\begin{corollary}\label{C-Randers-A0}
Let $\calf$ be a codimension-one totally umbilical $($for the metric $a)$ foliation
of a closed Randers space $(M^{m+1},\alpha+\beta)$ with constant sectional curvature~$\bar K$ of $a$.
If a nonzero vector field $\beta^\sharp\in\Gamma(T\calf)$ obeys $\bar\nabla\beta^\sharp=0$
then $\bar K=0$,
$\calf$ is totally geodesic and for $1\le k\le m$ $($for $k=1$, see also Remark~\ref{Ex-k1}$)$ we have
\begin{eqnarray}\label{E-IF-Randers-k-tg}
\nonumber
 && \int_M \Big(c^k\sigma_{k}(C^\sharp_\nu) -\frac12\,c^{-1}\,\<T_{k-1}(c\,C^\sharp_\nu)(\beta^\sharp),\, \bar Z\> \\
 && -\,\frac c2\,\big\<T_{k-1}(c\,C^\sharp_\nu
 -\frac12\,c^{-1}\bar Z^\flat\otimes\beta^\sharp)(\bar Z),\,\beta^\sharp\big\>
 \Big)\,{\rm d}V_a = 0.
\end{eqnarray}
\end{corollary}

\proof Since $\<\bar A(\beta^\sharp),\beta^\sharp\>=0$ (see the proof of Theorem~\ref{T-Randers-Kconst}),
we obtain $\bar H=0$.
Thus, \eqref{E-IF-Randers-k-tg} follows from \eqref{E-IF-Randers-k} with $\bar A=0$ and $\beta(\bar Z)=0$.
\qed

\section{Appendix: Invariants of a set of matrices}
\label{sec:inv}

Here, we collect the properties of the invariants
$\sigma_{{\lambda}} (A_1, \ldots , A_k)$ of real matrices $A_i$ that generalize the elementary symmetric functions
of a single symmetric matrix~$A$.
Let $S_k$ be the group of all permutations of $k$ elements.
 Given arbitrary quadratic $m\times m$ real matrices $A_1, \ldots  A_k$
and the unit matrix $I_{m}$, one can consider the determinant $\det(I_{m}+t_1A_1+\ldots+t_kA_k)$ and express it as a polynomial of
real variables ${\bf t}=(t_1, \dots  t_k)$. Given
 ${\lambda} = (\lambda_1, \ldots  \lambda_k)$, a sequence of nonnegative integers with
$|{\lambda}| := \lambda_1 + \ldots + \lambda_k\le m$, we shall
denote by $\sigma_{{\lambda}} (A_1, \ldots , A_k)$ its coefficient at
 ${\bf t}^{{\lambda}}=t_1^{\lambda_1}\cdot\ldots t_k^{\lambda_k}$:
\begin{equation}\label{eq00}
 \det(I_{m}+t_1A_1+\ldots+t_kA_k)=\sum\nolimits_{\,|{\lambda}|\,\le m}
 \sigma_{{\lambda}}(A_1, \ldots  A_k)\,{\bf t}^{{\lambda}}.
\end{equation}
Evidently, the quantities $\sigma_{{\lambda}}$ are invariants of conjugation by $GL(m)$-matrices:
\begin{equation}\label{eq01}
 \sigma_{{\lambda}}(A_1,\ldots A_k)=\sigma_{{\lambda}}(QA_1Q^{-1},\ldots QA_kQ^{-1})
\end{equation}
for all $A_i$'s, ${{\lambda}}$'s and nonsingular $m\times m$ matrices $Q$.
Certainly, $\sigma_i(A)$ (for a single symmetric matrix $A$)
coincides with the $i$-th elementary symmetric polynomial of the
eigenvalues $\{k_j\}$ of~$A$.

In the next lemma, we collect properties of these invariants.

\begin{lemma}[see \cite{rw2}]\label{lem1}
 For any ${\lambda}=(\lambda_1,\ldots \lambda_k)$ and any $m\times m$ matrices $A_i, A$  and $B$ one has

\noindent\
 (I) $\sigma_{{\lambda}}(0, A_2, \dots  A_{k})=0$ if $\lambda_1>0$
 and
 $\sigma_{0,\hat{{\lambda}}}(A_1, \dots  A_{k})
 =\sigma_{\hat{{\lambda}}}(A_2,\ldots A_{k})$
 where $\hat{{\lambda}}=(\lambda_2,\ldots \lambda_k)$,

\noindent\
 (II) $\sigma_{{\lambda}} (A_{s(1)},\ldots  A_{s(k)})=\sigma_{{\lambda} \circ s}
(A_1,\ldots A_k)$, where $s\in S_k$ and ${\lambda}\circ s=(\lambda_{s(1)},\ldots \lambda_{s(k)})$,

\noindent\
 (III) $\sigma_{{\lambda}}(I_{m},A_2,\ldots  A_k)
 =\genfrac{(}{)}{0pt}{}{m-|\hat{{\lambda}}|}{\lambda_1}\,
 \sigma_{\hat{{\lambda}}} (A_2,\ldots  A_k)$,

\noindent\
 (IV) $\sigma_{\lambda_1,\lambda_2,\,\hat{{\lambda}}}(A,A,A_3,\ldots  A_k)
 =\genfrac{(}{)}{0pt}{}{\lambda_1+\lambda_2}{\lambda_1}\,
 \,{\sigma_{\lambda_1+\lambda_2,\hat{{\lambda}}}}(A,A_3,\ldots  A_k)$,

\noindent\
 (V) $\sigma_{1,\hat{{\lambda}}}(A+B,A_2,\dots  A_{k})
 =\sigma_{1,\hat{{\lambda}}}(A,A_2,\ldots
A_{k})+\sigma_{1,\hat{{\lambda}}}(B,A_2,\ldots A_{k})$ and

 \quad
 $\sigma_{{\lambda}}(a\/A_1,A_2,\dots  A_{k})
 =a^{\lambda_1}\sigma_{{\lambda}}(A_1,A_2,\ldots A_{k})$ if $a\in\RR\setminus\{0\}$.
\end{lemma}

The invariants defined above can be used in calculation of the
determinant of a matrix $B(t)$ expressed as a power series
 $B(t)=\sum_{i=0}^\infty t^iB_i$.
Indeed, if one wants to express $\det(B(t))$ as a power series in
$t$, then the coefficient at $t^j$ depends only on the part
$\sum\nolimits_{\,i\le j} t^iB_i$ of $B(t)$.

\begin{lemma}[\cite{rw2}]\label{lem2}
 If $B(t)$, $t\in\RR$, is the $m\times m$ matrix given by
 $B(t)=\sum\nolimits_{i=0}^\infty t^iB_i$, $B_0=I_m$ then
\begin{equation}\label{eq03}
 \det(B(t))=1+\sum\nolimits_{k=1}^\infty\big(\sum\nolimits_{{\lambda},\|{\lambda}\|=k}
 \sigma_{{\lambda}}(B_1, \ldots  B_k)\big)\,{t}^k,
\end{equation}
where $\|{\lambda}\|=\lambda_1 + 2\lambda_2 + \ldots + k\lambda_k$
for ${\lambda}=(\lambda_1,\ldots \lambda_k)$. \hfill\qed
\end{lemma}

Since $\det:{\cal M}(m)\to\RR,\ {\cal M}(m)\approx\RR^{m^2}$ being
the space of all $m\times m$-matrices, is a polynomial function, the
series in (\ref{eq03}) is convergent for all $t\in (-r_0, r_0)$,
where $r_0=1/\limsup\nolimits_{\,k\to\infty}\| B_k\|^{1/k}$ is the
radius of convergence of the series $B(t)$.

By the First Fundamental Theorem of Matrix Invariants, see \cite{gu}, all the invariants $\sigma_{{\lambda}}$
can be expressed in terms of the traces of the matrices involved and their products.

\begin{lemma}[\cite{rw2}]
 For arbitrary matrices $B$, $C$ and $k,l>0$ we have
\begin{equation*}
 \sigma_{k,l}(B,C)=\sigma_k(B)\,\sigma_l(C) -\sum\nolimits_{\,i=1}^{\,\min(k,l)}\sigma_{k-i,l-i,i}(B,C,BC).
\end{equation*}
In particular, for $l=1$, it follows that
\begin{eqnarray}\label{eq-sigma-k1}
 \sigma_{k,1}(B,C) \eq \sum\nolimits_{i=0}^k(-1)^i\sigma_{k-i}(B)\,{\tr}(B^iC) =\tr(T_{k}(B)C).
\end{eqnarray}
\end{lemma}

\begin{lemma}\label{L-AplusB}
Let $A,C$ be $m\times m$ matrices and ${\rm rank}\, A=1$. Then
\begin{equation}\label{E-CArank1}
 \sigma_k(C+A) = \sigma_k(C) +\tr(T_{k-1}(C)A).
\end{equation}
\end{lemma}

\proof
There exists a nonsingular matrix $Q$ such that $\tilde A=QAQ^{-1}$ has one nonzero element, $\tilde a_{1i}\ne0$ for some $i$
(the simplest rank one matrix).
By \eqref{eq01}, $\sigma_{k,l}(\tilde C,\tilde A)=\sigma_{k,l}(C, A)$ where $\tilde C=QCQ^{-1}$.
By Laplace's formula (which expresses the determinant of a matrix in terms of its minors),
$\det(I_m+t\tilde C+s\tilde A)$ is a linear function in $s\in\RR$; hence, see \eqref{eq00}, $\sigma_{k,l}(\tilde C,\tilde A)=0$ for $l>1$.
By the above, $\sigma_{k,l}(C, A)=0$ for $l>1$ and all $k$.
 Using the identity, see \cite{rw2},
\begin{equation}\label{E-sigma22}
 \sigma_k(C_1+C_2)=\sum\nolimits_{\,i=0}^{\,k}\sigma_{\,k-i, i}(C_1,C_2),
\end{equation}
we find that
\begin{eqnarray*}
 \sigma_k(C+A) \eq \sigma_k(C) +\sigma_{k-1,1}(C, A).
\end{eqnarray*}
By \eqref{eq-sigma-k1}, $\sigma_{k-1,1}(C, A)=\tr(T_{k-1}(C)A)$ and \eqref{E-CArank1} follows.
\qed

\begin{corollary}\label{C-AplusB}
Let $C,D,A_i$ be $m\times m$ matrices and ${\rm rank}\, A_i=1\ (1\le i\le s)$. Then
\begin{eqnarray}\label{E-sigma-k}
\nonumber
 && \sigma_k(C+D+A_1+\ldots A_s) = \sigma_k(C) +\sum\nolimits_{j>0}\sigma_{k-j,j}(C, D) \\
 \plus\tr(T_{k-1}(C+D)A_1) +\ldots+\tr(T_{k-1}(C+D+A_1+\ldots+A_{s-1})A_s).
\end{eqnarray}
\end{corollary}

\proof This follows from Lemma~\ref{L-AplusB} and \eqref{eq-sigma-k1}. For $s=1$, we obtain
\begin{eqnarray*}
 &&\hskip1mm\sigma_k(C+D+A_1) \overset{\eqref{E-CArank1}}{=} \sigma_k(C+D) +\tr(T_{k-1}(C+D)A_1) \\
 && \overset{\eqref{E-sigma22}}{=} \sigma_k(C) +\sum\nolimits_{j>0}\sigma_{k-j,j}(C, D) +\tr(T_{k-1}(C+D)A_1).
\end{eqnarray*}
Then, by induction for $s$, \eqref{E-sigma-k} follows.
\qed

\smallskip

Let $C_i$ and $P_i$ be $m$-vectors (columns) and $I_m$ the identity $m$-matrix and $1\le i\le j\le m$.
Note that $C_i P_j^t$ are $m\times m$-matrices of rank 1 with
\begin{eqnarray*}
 \sigma_1(C_i P_j^t)=C_i^t P_j=P_j^t C_i,\quad
 \sigma_2(C_i P_j^t)=0,\\
 (I_m+C_i P_j^t)^{-1} = I_m - (1+C_i^t P_j)^{-1}\, C_i P_j^t.
\end{eqnarray*}

\begin{lemma}
We have
 $\det(I_m +\sum\nolimits_{i=1}^k C_i P_i^t) = 1 + \det(\{C_i^t P_j\}_{1\le i,j\le k})$.
For example,
\begin{eqnarray*}
 && \det(I_m+C_1 P_1^t) = 1+C_1^t P_1\,,\\
 && \det(I_m+C_1 P_1^t +C_2 P_2^t) = 1+C_1^t P_1+C_2^t P_2
 +C_1^t P_1\cdot C_2^t P_2 -C_1^t P_2\cdot C_2^t P_1\,,
\end{eqnarray*}
and so on.
\end{lemma}

%
%
%
%
%


\begin{thebibliography}{999.}%

\bibitem{arw2014}
 K. Andrzejewski, V. Rovenski and P. Walczak, {\em Integral formulas in foliation theory}, 73-82,
in ``\textit{Geometry and its Applications}", Springer Proc. in Math. and Statistics, 72, Springer,~2014.

\bibitem{bcs}
 D. Bao, S.S. Chern, Z. Shen, {\em An Introduction to Riemann-Finsler geometry}, Springer,~2000.

\bibitem{blr}
 F. Brito, R. Langevin,  H. Rosenberg, {\em Int\'{e}grales de courbure sur des vari\'{e}t\'{e}s feuillet\'{e}es}.
J.~Diff. Geom. {16} (1981), 19-50.

\bibitem{cs}
 X. Cheng,  Z. Shen, {\em Finsler geometry. An approach via Randers spaces}, Springer, 2012.

\bibitem{cs2}
 S.S. Chern, Z. Shen, {\em Riemann-Finsler geometry}, World Scientific, 2005.

\bibitem{gu}
 G.G. Gurevich, {\em Foundations of the theory of algebraic invariants}, Noordhof, Groningen, 1964.

\bibitem{lw}
R. Langevin, P. Walczak, {\em Conformal geometry of foliations}, Geom. Dedic.  {132} (2008), 135-178.

\bibitem{mo}
 X. Mo, {\em An introduction to Finsler geometry}, World Scientific, 2006.

\bibitem{ra}
 G. Randers, {\em On an asymmetrical metric in the four-space of general relativity}, Phys. Rev. 59, (1941) 195-199.

\bibitem{re}
 G. Reeb, {\em Sur la courbure moyenne des vari\'{e}t\'{e}s int\'{e}grales d'une \'{e}quation de Pfaff} $\omega=0$.
 C. R. Acad. Sci. Paris, {231}, (1950) 101-102.

\bibitem{ro}
 V. Rovenski, {\em Foliations on Riemannian manifolds and submanifolds}, Birkh\"{a}user, Boston 1998.

\bibitem{rw1}
 V. Rovenski, P. Walczak, {\em Topics in extrinsic geometry of codimension-one foliations}, Springer 2011.

\bibitem{rw2}
 V. Rovenski, P. Walczak, {\em Integral formulae on foliated symmetric spaces}, Math. Ann. {352}, (2012) 223-237.

\bibitem{sh1}
 Z. Shen, {\em On Finsler geometry of submanifolds}, Math. Ann. {311}, (1998) 549-576.

\bibitem{sh2}
 Z. Shen, {\em Lectures on Finsler geometry}, World Scientific Publishers, 2001.

\bibitem{topp}
P. Topping, {\em Lectures on the Ricci flow}, LMS Lecture Notes 325, Cambridge Univ. Press, 2006.

\end{thebibliography}
\end{document}